\documentclass[11pt,english,oneside]{amsart}
\usepackage[T1]{fontenc}
\usepackage[latin9]{inputenc}
\usepackage{geometry}
\usepackage{amssymb}
\usepackage{color}

\usepackage{graphicx,color}
\usepackage{amsmath, amssymb, graphics}

\usepackage{graphicx}

\makeatletter
\numberwithin{equation}{section} 
\numberwithin{figure}{section} 
  \@ifundefined{theoremstyle}{\usepackage{amsthm}}{}
  \theoremstyle{plain}
  \newtheorem{thm}{Theorem}[section]
  \theoremstyle{plain}
  \newtheorem{cor}[thm]{Corollary}
  \theoremstyle{plain}
  
  \theoremstyle{remark}
  
  \theoremstyle{remark}
  
  \theoremstyle{plain}
  \newtheorem{lem}[thm]{Lemma}

\usepackage{geometry}

\def\bfR#1{{\bf R}^#1}

\def\com#1{ \hbox{#1}}

\smallskip
\def\<{{\langle }}
\def\>{{\rangle }}

\makeatother

\usepackage{babel}

\def\bfR#1{{\bf R}^#1}

\def\com#1{ \quad\hbox{#1}\quad}

\smallskip
\def\<{{\langle }}
\def\>{{\rangle }}

\makeatother

\begin{document}

\title{Embedded constant mean curvature hypersurfaces on spheres}

\author{ Oscar M. Perdomo }

\date{\today}

\curraddr{Department of Mathematics\\
Central Connecticut State University\\
New Britain, CT 06050\\
}

\email{ perdomoosm@ccsu.edu}

\begin{abstract}

Let $m\ge2$ and $n\ge2$ be any pair of integers. In this paper we prove that if $H$ is between the numbers $ \cot(\frac{\pi}{m})$ and $b_{m,n}=\frac{(m^2-2)\sqrt{n-1}}{n\sqrt{m^2-1}}$,
then, there exists a non isoparametric, compact embedded hypersurface in $S^{n+1}$ with constant mean curvature $H$ that admits the group $O(n)\times Z_m$ in their group of isometries, here $O(n)$ is the set of $n\times n$ orthogonal matrices and $Z_m$ are the integers mod $m$. When $m=2$ and $H$ is close to the boundary value $0= \cot(\frac{\pi}{2})$, the hypersurfaces look like two very close  $n$-dimensional spheres with two catenoid necks attached, similar to  constructions made by Kapouleas. When $m>2$ and $H$ is close to $\cot(\frac{\pi}{m})$, the hypersurfaces look like a necklets made out of $m$ spheres with $m+1$ catenoid necks attached, similar to constructions made by Butscher and Pacard. In general, when $H$ is close to $b_{m,n}$ the hypersurface is close to an isoparametric hypersurface with the same mean curvature. As a consequence of the expression of these bounds for $H$, we have that every $H$ different from $0,\pm\frac{1}{\sqrt{3}}$ can be realized as the mean curvature of a non isoparametric CMC surface in $S^3$. For hyperbolic spaces we prove that every non negative $H$ can be realized as the mean curvature of an embedded CMC hypersurface in $H^{n+1}$, moreover we prove that when $H>1$ this hypersurface admits the group $O(n)\times Z$ in its group of isometries. Here $Z$ are the integer numbers. As a corollary of the properties proven for these hypersurfaces, for any $n\ge 6$, we construct non isoparametric compact minimal hypersurfaces in $S^{n+1}$ which cone in $\bfR{{n+2}}$ is stable. Also, we will prove that the stability index of every non isoparametric minimal hypersurface  with  two principal curvatures in $S^{n+1}$ is greater than $2n+5$.

\end{abstract}


\maketitle
\section{Introduction}

Minimal hypersurfaces on spheres with exactly two principal curvatures everywhere were studied by Otsuki in (\cite{O}), he reduced the problem of classifying them all, to the problem of solving an ODE, and the problem of deciding about their compactness, to the problem of studying a real function given in term of an integral that related two periods of two functions involved in the  immersions that he found. In this paper we changed the minimality condition for the constant mean curvature condition and using a slightly different point of view, we got similar results. As pointed out by Do Carmo and Dajczer in (\cite{D-D}), these hypersurfaces are rotations of a plane profile curve.  The existence of these  hypersurfaces as immersions have been established  in (\cite{D-D}) and in (\cite{SI}). Partial result about the condition for small values of $H$ that guarantee embedding were found in (\cite{B-L}). The main work in this paper consists in studying these profile curves and in deciding when they are embedded. Lemma (\ref{the lemma}) and its corollary (\ref{cor of the lemma})  played an important role in the understanding of the period of these profile curves for they were responsible of getting explicit formulas for these two numbers

$$a_{n,m}=\cot{\frac{\pi}{m}} \com{and} b_{n,m}=\frac{ (m^2-2)\, \sqrt{n-1}}{n \sqrt{m^2-1}} $$

with the property that every $H$ between  $a_{n,m}$ and $b_{n,m} $ can be realized as the mean curvature of a non isoparametric embedded CMC hypersurface in $S^{n+1}$, such that its profile curve is invariant under the group of rotation by an angle $\frac{2\pi}{m}$, and therefore the hypersurface admits the group $O(n)\times Z_m$ in its group of isometries.

 One of the differences between the analysis of the profile curve in this work and Otsuki's, is that Otsuki used the supporting function of this profile curve. We, instead, studied the radius and the angle separately, and we proved that the angle function is increasing, which helps us to decide when this curve is injective and consequently, when the immersion is an embedding. In order to understand this angle function we got some help from the understanding of three vector fields defined in section (\ref{section vector fields that lie on a planes}).

Since we can obtain a similar formula for the angle of the profile curve for CMC hypersurfaces in the hyperbolic space, we also extend the result in this case in order to explicitly exhibit embedded examples the hyperbolic space. Also, similar results are obtained in Euclidean spaces. These results on embedded hypersurfaces on hyperbolic spaces and Euclidean spaces where proven in (\cite{H}) with different techniques.

As a consequence of the symmetries proven for all compact constant mean curvatures in $S^{n+1}$ with two principal curvatures everywhere, we proved that all these examples with $H=0$, have stability index greater than $2n+5$. In this direction there is a conjecture that  states that the only minimal hypersurfaces in $S^{n+1}$ with stability index $2n+5$ are the isoparametric with two principal curvatures. Some partial results for this conjecture were proven in (\cite{P1}). Also, since it is not difficult to prove that these examples can be chosen to be as close as we want from the isoparametric examples, we proved that some of Otsuki's minimal hypersurfaces are examples of non isoparametric compact stable minimal truncated cones in $\bfR{{n+2}}$ for $n\ge6$.

The author would like to express his gratitude to Professor Bruce Solomon for discussing the hypersurfaces with him and pointing out the similarity between them and the Delaunay's surfaces.

\section{Preliminaries}

\medskip

Let $M$ be an n-dimensional hypersurface of the $(n+1)$-dimensional unit sphere $S^{{n+1}}\subset \bfR{{n+2}}$. Let $\nu:M\to S^{n+1}$ be a Gauss map and $A_p:T_pM\to T_pM$ the shape operator, notice that

$$A_p(v)=-\bar{\nabla}_v\nu\com{for all} v\in T_pM $$

where $\bar{\nabla}$ is the Euclidean connection in $\bfR{{n+2}}$. We will denote by $||A||^2$ the square of the norm of the shape operator.

If $X$, $Y$ and $Z$ are vector fields on $M$, $\nabla_XY$  represents the Levi-Civita connection on $M$ with respect to the metric induced by $S^{n+1}$ and $[X,Y]=\nabla_XY-\nabla_YX$ represents the Lie bracket, then, the curvature tensor on $M$ is defined by

\begin{eqnarray}
\label{definition of curvature tensor}
R(X,Y)Z=\nabla_Y\nabla_XZ-\nabla_X\nabla_YZ+\nabla_{[X,Y]}Z
\end{eqnarray}

and the covariant derivative of $A$ is defined by
\begin{eqnarray}
\label{definition of DA}
DA(X,Y,Z)=Z\<A(X),Y\>-\<A(\nabla_ZX),Y\>-\<A(X),\nabla_ZY\>
\end{eqnarray}

the Gauss equation is given by,

\begin{eqnarray}
\label{Gauss equation}
R(X,Y)Z=\<X,Z\>Y-\<Y,Z\>X+\<A(X),Z\>A(Y)-\<A(Y),Z\>X
\end{eqnarray}

and the Codazzi equations are given by,

\begin{eqnarray}
\label{Codazzi equation}
DA(X,Y,Z)=DA(Z,Y,X)
\end{eqnarray}

Let us denote by $\kappa_1,\dots,\kappa_n$ the principal curvatures of $M$ and, by $\displaystyle{H=\frac{\kappa_1+\dots+\kappa_n}{n}}$ the mean curvature of $M$. We will assume that $M$ has exactly two principal curvatures everywhere and that $H$ is a constant function on $M$. More precisely, we will assume that

$$\kappa_1=\dots=\kappa_{n-1}=\lambda,\quad\kappa_n=\mu \com{and} (n-1)\lambda+\mu=nH$$

By changing $\nu$ by $-\nu$ if necessary we can assume without loss of generality that $\lambda-\mu> 0$. Let $\{e_1,\dots,e_n\} $ denotes a locally defined orthonormal frame such that

\begin{eqnarray}
\label{frame}
 A(e_i)=\lambda e_i\com{for $i=1,\dots,n-1$ }
\com{and}A(e_n)=\mu e_n
\end{eqnarray}

The next Theorem is well known, see (\cite{O}), for completeness sake and as part of preparation for the deduction of other formulas, we will show a proof here,

\begin{thm}
\label{equations}

 If $M\subset S^{n+1}$ is a CMC hypersurface with two principal curvatures and dimension greater than 2, and  $\{e_1,\dots,e_n\} $ is a locally defined orthonormal frame such that (\ref{frame}) holds true, then,

 \begin{eqnarray*}
v(\lambda)&=&0 \com{for any } v\in \hbox{Span}\{e_1,\dots,e_{n-1}\}\\
\nabla_ve_n&=&\frac{e_n(\lambda)}{\mu-\lambda}\, v\com{for any } v\in \hbox{Span}\{e_1,\dots,e_{n-1}\}\\
\nabla_{e_n}e_n&=&0\\
1+\lambda\mu&=&e_n(\frac{e_n(\lambda)}{\lambda-\mu})-(\frac{e_n(\lambda)}{\lambda-\mu})^2
 \end{eqnarray*}
$$[e_i,e_j] \in  \hbox{Span}\{e_1,\dots,e_{n-1}\}  \com{for any } i,j\in\{1,\dots,n-1\}$$

\end{thm}
\begin{proof}
For any $i,j\in\{1,\dots,n-1\}$ with $i\ne j$ (here we are using the fact that the dimension of $M$ is greater than 2) and any $k\in\{1,\dots,n\}$, we have that,

 \begin{eqnarray*}
 DA(e_i,e_j,e_k)&=&e_k\<\, A(e_i),e_j\,\> -\<\, A(\nabla_{e_k}e_i),e_j\,\> -\<\, A(e_i),\nabla_{e_k}e_j\,\> \\
 &=&e_k(\lambda \<\, e_i,e_j\,\> )-\<\, \nabla_{e_k}e_i,A(e_j)\,\> -\lambda\<\, e_i,\nabla_{e_k}e_j\,\> \\
&=& e_k( 0)-\lambda\<\, \nabla_{e_k}e_i,e_j\,\> -\lambda\<\, e_i,\nabla_{e_k}e_j\,\> \\
&=& 0-\lambda e_k(\<\, e_i,e_j\,\> )\\
&=& 0
\end{eqnarray*}

On the other hand,
\begin{eqnarray*}
DA(e_i,e_i,e_j)&=&e_j\<\, A(e_i),e_i\,\> -\<\, A(\nabla_{e_j}e_i),e_i\,\> -\<\, A(e_i),\nabla_{e_j}e_i\,\> \\
 &=&e_j(\lambda)-\lambda e_j(\<\, e_i,e_i\,\> )\\
&=& e_j(\lambda)
\end{eqnarray*}

By Codazzi equation (\ref{Codazzi equation}), we get that $e_j(\lambda)=0$, for all $j\in\{1,\dots,n-1\}$, therefore $v(\lambda)=0$ for any  $v\in \hbox{Span}\{e_1,\dots,e_{n-1}\}$. Now,
\begin{eqnarray*}
DA(e_i,e_n,e_j)&=&e_j\<\, A(e_i),e_n\,\> -\<\, A(\nabla_{e_j}e_i),e_n\,\> -\<\, A(e_i),\nabla_{e_j}e_n\,\> \\
&=&e_j(\lambda\<\, e_i,e_n\,\> )-\<\, \nabla_{e_j}e_i,A(e_n)\,\> -\lambda\<\, e_i,\nabla_{e_j}e_n\,\> \\
&=&e_j(0)-\mu\<\, \nabla_{e_j}e_i,e_n\,\> -\lambda\<\, e_i,\nabla_{e_j}e_n\,\> +( \lambda\<\, \nabla_{e_j}e_i,e_n\,\> -\lambda\<\, \nabla_{e_j}e_i,e_n\,\>)\\
 &=&(\lambda-\mu)\<\, \nabla_{e_j}e_i,e_n\,\> -\lambda\ e_j(\<\, e_i,e_n\,\> )\\
&=& (\mu-\lambda)\<\, e_i,\nabla_{e_j}e_n\,\>
\end{eqnarray*}

Since $\mu-\lambda>0$, using Codazzi equations  we get that,

\begin{eqnarray}\label{den1}
\<\, e_i,\nabla_{e_j}e_n\,\> =0\com{for any $i,j\in\{1,\dots,n-1\}$ with $i\ne j$}
\end{eqnarray}

Now, for any $i\in\{1,\dots,n-1\}$, using the same type of computations as above we can prove that,

$$DA(e_i,e_i,e_n)=e_n(\lambda)\, =\, DA(e_i,e_n,e_i)=(\mu-\lambda)\<\, e_i,\nabla_{e_i}e_n\, \>$$

and

$$DA(e_n,e_n,e_i)=e_i(\mu)\, = \, 0\,=\, DA(e_i,e_n,e_n)=(\mu-\lambda)\<\, e_i,\nabla_{e_n}e_n\, \>$$

Therefore,

\begin{eqnarray}\label{den2}
\<\, e_i,\nabla_{e_i}e_n\,\> =\frac{e_n(\lambda)}{\mu-\lambda}\com{and}\<\, e_i,\nabla_{e_n}e_n\,\>=0\com{for any $i\in\{1,\dots,n-1\}$ }
\end{eqnarray}

Since $e_n$ is a unit vector field, we have that $\<\,\nabla_{e_k}e_n,e_n\,\>=0$ for any $k$. From the equations (\ref{den1} ) and (\ref{den2} ) we conclude that

$$\displaystyle{\nabla_ve_n=\frac{e_n(\lambda)}{\mu-\lambda}}\, v\com{for any}v\in \hbox{Span}\{e_1,\dots,e_{n-1}\}\com{and} \nabla_{e_n}e_n=0$$

Notice that for any $i,j\in\{1,\dots,n-1\}$ with $i\ne j$, using equation (\ref{den1}), we get that

$$\<\, [e_i,e_j],e_n\,\> =\<\, \nabla_{e_i}e_j-\nabla_{e_j}e_i,e_n\,\>=\<\,e_i,\nabla_{e_j}e_n\,\>-\<\,e_j,\nabla_{e_i}e_n\,\>=0$$

Therefore $[e_i,e_j]\in\hbox{Span}\{e_1,\dots,e_{n-1}\}$. Finally we will use Gauss equation to prove the differential equation on $\lambda$. First let us point out that, using the equation (\ref{den2}), we can prove that $\<\, [e_n,e_1],e_n\,\> =0$ and therefore we have that $[e_n,e_1]\in\hbox{Span}\{e_1,\dots,e_{n-1}\}$. Now, by Gauss equation we get,

\begin{eqnarray*}
1+\lambda\mu  &=&\<\, R(e_n,e_1)e_n,e_1\,\> \\
&=&\<\,\nabla_{e_1}\nabla_{e_n}e_n-\nabla_{e_n}\nabla_{e_1}e_n+\nabla_{[e_n,e_1]}e_n,e_1\, \>\\
&=&\<\,0-\nabla_{e_n}(\frac{e_n(\lambda)}{\mu-\lambda}\, e_1)+ \frac{e_n(\lambda)}{\mu-\lambda}\, [e_n,e_1]   \, ,e_1 \>\\
&=& -e_n(\frac{e_n(\lambda)}{\mu-\lambda})+\frac{e_n(\lambda)}{\mu-\lambda} \<\,\nabla_{e_n}e_1-\nabla_{e_1}e_n ,e_1  \, \>\\
&=&-e_n(\frac{e_n(\lambda)}{\mu-\lambda})-(\frac{e_n(\lambda)}{\mu-\lambda})^2\\
&=&e_n(\frac{e_n(\lambda)}{\lambda-\mu})-(\frac{e_n(\lambda)}{\lambda-\mu})^2
\end{eqnarray*}

\end{proof}

\section{Construction of the examples}
 We will maintain the notation of the previous section and we will prove a serious of identities and results that will allow us to easier state and prove the theorem that defines the examples at the end of this section.

\subsection{The function $w$ and its solution along a line of curvature\label{section the function w and its solutions}}

Since $(n-1)\lambda+\mu=n H$, we get that

\begin{eqnarray}
\label{definition of w}
\lambda-\mu=\lambda-(nH-(n-1)\lambda)=n(\lambda-H)=nw^{-n}\com{where} w=(\lambda-H)^{-\frac{1}{n}}
\end{eqnarray}

Recall that we are assuming that $\lambda-\mu$ is always positive, then, $w$ is a smooth differentiable function. By the definition of $w$ given in (\ref{definition of w}) we have that,

\begin{eqnarray}
\label{derivative of w}
e_n(w)=-\frac{1}{n}(\lambda-H)^{-\frac{n+1}{n}}e_n(\lambda)=-\frac{1}{n}w^{n+1}\, e_n(\lambda)=-w \frac{e_n(\lambda)}{\lambda-\mu}
\end{eqnarray}

The second order differential equation in Theorem (\ref{equations}) can be written using the function $w$ as,

\begin{eqnarray}
\label{second order equation of w first version}
e_n\big{(}\,\frac{e_n(w)}{w}\,\big{)}+\big{(}\,\frac{e_n(w)}{w}\,\big{)}^2+1+\lambda \mu=0
\end{eqnarray}

and if we write $\lambda $ and $\mu$ in terms of $w$ we get

\begin{eqnarray}
\label{second order equation of w second version bis}
e_n(\big{(}\,\frac{e_n(w)}{w}\,\big{)})+\big{(}\,\frac{e_n(w)}{w}\,\big{)}^2-
\frac{(n-1)}{w^{2n}}-\frac{(n-2)H }{w^n}+H^2+1=0
\end{eqnarray}

Deriving the previous equation, we have used the following identities,
\begin{eqnarray}
\label{lambda in terms of w}
\lambda=w^{-n}+H\com{and} \mu=H-(n-1) w^{-n}
\end{eqnarray}

From the Equation (\ref{derivative of w})  we get that
\begin{eqnarray}
\label{derivative of lambda}
e_n(\lambda)=-(\lambda-\mu)\, \frac{e_n(w)}{w}
\end{eqnarray}

The equation above allows us to write one of the equations in Theorem (\ref{equations}) as
\begin{eqnarray}
\label{derivative of lambda e_n in the direction e_i}
\bar{\nabla}_{v}e_n= \frac{e_n(w)}{w} \, v\com{for any } v\in \hbox{Span}\{e_1,\dots, e_{n-1}\}
\end{eqnarray}

Notice that equation  (\ref{second order equation of w second version bis}) reduces to

\begin{eqnarray}
\label{second order equation of w second version}
\frac{e_n(e_n(w))}{w}-\frac{(n-1)}{w^{2n}}-\frac{(n-2)H }{w^n}+H^2+1=0
\end{eqnarray}

and therefore multiplying by $2w e_n(w)$ we have that there exists a  constant $C$
such that,

\begin{eqnarray}
\label{integral equation for w} (e_n(w))^2+w^{2-2n}+(1+H^2)w^2+2Hw^{2-n}=C
\end{eqnarray}

Let us denote by $x:M\to \bfR{{n+2}}$ the position viewed as a map, and by $\bar{\nabla}$ the Euclidean connection on $\bfR{{n+1}}$. Using the equations in Theorem (\ref{equations}) and the fact that $\bar{\nabla}_{v}x=v$,  $\<\,x,\nu(x)\,\>=0$ and $\<\,\nu(x),\nu(x)\,\>=1$, we get that
\begin{eqnarray}
\label{derivatives in the Euclidean space 1}
\bar{\nabla}_{e_n}e_n&=&-x+\mu \nu\\
\label{derivatives in the Euclidean space 2}\bar{\nabla}_{e_n}\nu &=& -\mu e_n                   \\
\label{derivatives in the Euclidean space 3}\bar{\nabla}_{e_n}x   &=&   e_n
\end{eqnarray}

Let us fix a point $p_0\in M$, and let us denote by $\gamma(u)$ the only geodesic in $M$ such that $\gamma(0)=p_0$ and $\gamma^{\prime}(0)=e_n(p_0)$. Since $\nabla_{e_n} e_n$ vanishes, then $\gamma(u)=e_n(\gamma(u))$. Notice that $\gamma(u)$ is also a line of curvature. Let us denote by $g(u)=w(\gamma(u))$. Equation (\ref{integral equation for w}) implies that

\begin{eqnarray}
\label{integral equation for w 1} (g^\prime)^2+g^{2-2n}+(1+H^2)g^2+2Hg^{2-n}=C
\end{eqnarray}

 or equivalently,

\begin{eqnarray}
\label{integral equation for w 2} \frac{g^{n-1}\, g^\prime}{\sqrt{C g^{2n-2}-1-(1+H^2)g^{2n}-2Hg^{n}}}=\pm 1
\end{eqnarray}

 It is clear that the constant $C$ must be positive and moreover, in order to solve this equation we need to consider a constant $C$ such that the polynomial

 \begin{eqnarray}\label{the polynomial}
  \xi(s)=C s^{2n-2}-1-(1+H^2)s^{2n}-2Hs^{n}
 \end{eqnarray}
is positive on a interval $(t_1,t_2)$ with $0<t_1<t_2$ and $\xi(t_1)=0=\xi(t_2)$. Notice that for every $H$ it is possible to pick $C$ such that $\xi$ is positive on an interval because $\xi$ is a polynomial of even degree with negative leading coefficient, $\xi(0)=-1$,  and if $C$ is big enough, this polynomial takes positive values for positive values of $s$.
Let us assume that $t_1$ and $t_2$ are as above and moreover let us assume that $\xi^{\prime}(t_1)$ and $\xi^{\prime}(t_2)$ are not zero, so that the function $G$ that we are about to define is well defined on $[t_1,t_2]$. Let us consider

$$G(s)=\int_{t_1}^s\frac{t^{n-1}}{\sqrt{C t^{2n-2}-1-(1+H^2)t^{2n}-2Ht^{n}}}dt\com{for} t_1\le u\le t_2$$

and let $T=2G(t_2)$. Since $G^\prime(s)>0$ for $s\in(t_1,t_2)$, then we can consider the inverse of the function $G$. Denoting by $F:[0,\frac{T}{2}]\to [t_1,t_2]$ the inverse of $G$, a direct verification shows that the $T$-periodic  function that satisfies

$$g(u)=F(u)\com{for}0\le u\le \frac{T}{2}\com{and} g(u)=F(T-u)\com{for} \frac{T}{2}\le u\le T$$
is a solution of the equation (\ref{integral equation for w 2})

\subsection{The vector field $\eta$ \label{properties of eta}}

Let us define the following normal vector field along $M$

 $$\eta=-\frac{e_n(w)}{w}\, e_n+\lambda\, \nu-x$$

 The vector field $\eta$ has the following properties

 \begin{enumerate}
 \item $\displaystyle{\<\eta,\eta\>=(\frac{e_n(w)}{w})^2+\lambda^2+1=\frac{C}{w^2}}$, this follows from  Equation (\ref{integral equation for w}) and the definition of $\lambda$ in terms of $w$, (\ref{lambda in terms of w}).
\item  $\displaystyle{\bar{\nabla}_{e_n}\eta=-\frac{e_n(w)}{w} \eta}$. This is the equation that makes all the constructions work in this section. The equation follows from Equations (\ref{derivatives in the Euclidean space 1}), (\ref{derivatives in the Euclidean space 2}) and (\ref{derivatives in the Euclidean space 3} ) and the first and second differential equations for the function $w$, especially, Equation (\ref{second order equation of w first version}) and Equation (\ref{derivative of lambda}).
    \item For any $i\in\{1,\dots ,n-1\}$,  $\bar{\nabla}_{e_i}(x+\frac{w^2}{C}\, \eta)$ vanishes. The proof of this identity is similar, and additionally, uses the Equation (\ref{derivative of lambda e_n in the direction e_i}).
        \item $\displaystyle{\<x+\frac{w^2}{C}\, \eta,x+\frac{w^2}{C}\, \eta\>=1-\frac{w^2}{C}}$.
 \end{enumerate}

\subsection{Vector fields that lies on a plane \label{section vector fields that lie on a planes} }

Now that we have computed the function $g(u)=w(\gamma(u)) $, let us understand better the geodesic $\gamma$. The equations (\ref{derivatives in the Euclidean space 1}), (\ref{derivatives in the Euclidean space 2}) and (\ref{derivatives in the Euclidean space 3} ) imply that

$$X(u)=e_n(\gamma(u)),\quad Y(u)=\nu(\gamma(u))\com{and} Z(u)=\gamma(u) $$

satisfy an  ordinary linear differential equation in the variable $u$ with periodic coefficients (notice that $\mu(\gamma(u))$ is a function of $g(u)$). By the existence and uniqueness theorem of ordinary differential equations we get the solutions
$X(u)$, $Y(u)$ and $Z(u)$ lies in the three dimensional space
\begin{eqnarray}\label{the three dimensional space}
\Gamma_{p_0}=\hbox{Span}\{e_n(p_0),\nu(p_0),p_0\}
\end{eqnarray}

For the sake of simplification, we will consider the $T$-periodic function $r:{\bf R}\to {\bf R}$ defined by

$$ r(u)=\frac{g(u)}{\sqrt{C}}$$

It is not difficult to check that the function $r$ satisfies the following equations

\begin{eqnarray}\label{equation for r}
\frac{r^{\prime\prime}}{r}+1+\lambda\mu=0,\quad (r^\prime)^2+r^2\, (1+\lambda^2)=1,\quad
\lambda^\prime=-(\lambda-\mu)\frac{r^\prime}{r}
\end{eqnarray}

In the previous equations we are abusing of the notation with the name of the functions $\lambda$ and $\mu$, in this case, and whenever is understood by the context, they will also denote the function $\lambda(\gamma(u))$ and $\mu(\gamma(u))$ respectively. In the construction of the examples that we will be considering, the function $1-r^2$ needs to be positive. We can achieve this by assuming that $H\ge 0$ because this condition will imply that $\lambda>0$, and therefore $r<1$.

Let us define the following vector fields along $\gamma$

\begin{eqnarray*}
B_1(u)&=&\eta(\gamma(u))=-\frac{r^\prime}{r}\, X+ \lambda \, Y- Z\\
B_2(u)&=& -\frac{rr^\prime}{\sqrt{1-r^2}} \, X+\frac{r^2\lambda}{\sqrt{1-r^2}}\, Y+\sqrt{1-r^2} Z\\
B_3(u)&=& \frac{r \lambda}{\sqrt{1-r^2}}\, X+\frac{r^\prime}{\sqrt{1-r^2}}\, Y
\end{eqnarray*}

Using the equations in section (\ref{properties of eta}), the Equations (\ref{derivatives in the Euclidean space 1}), (\ref{derivatives in the Euclidean space 2}) and (\ref{derivatives in the Euclidean space 3} ) that give the derivative of  the vector fields $X$, $Y$ and $Z$,  and Equation (\ref{equation for r}),  we can check the following properties.

\begin{enumerate}
\item  $B_1(u)$, $B_2(u)$ and $B_3(u)$ lie on the three dimensional subspace $\Gamma_{p_0}$
\vskip.5cm
\item $${B_1^\prime=-\frac{r^\prime}{r} B_1}$$

\item $${\<B_1,B_2\>=0,\quad \<B_1,B_3\>=0\com{and}  \<B_2,B_3\>=0  }$$

 \item $${\<B_2,B_2\>=1,\quad \<B_3,B_3\>=1 \com{and} \<B_1,B_1\>=\frac{1}{r^2} }$$

 \item From the previous items we get that
 $$B_2^\prime= h B_3\com{ and} B_3^\prime=-hB_2\com{ for some function }h:{\bf R}\to {\bf R}$$
 These equations hold true because

      $$\<B_2^\prime,B_1\>=-\<B_1^\prime,B_2\>=\frac{r^\prime}{r}\<B_1,B_2\>=0\com{likewise} \<B_3^\prime,B_1\>=0$$

\item From the previous item we get  that the vectors $B_2$ and $B_3$ lie on a two dimensional subspace.

 \item We have that

$$\<B_3^\prime ,Z\>=-\frac{r\lambda}{\sqrt{1-r^2}}\com{and} \<B_2,Z\>=\sqrt{1-r^2}$$

therefore the function $h$ in the previous item is given by  $\frac{r\lambda}{1-r^2}$. It follows that,

$$B_2^\prime= \frac{r\lambda}{1-r^2} B_3\com{and} B_3^\prime=-\frac{r\lambda}{1-r^2}B_2$$

The fact that $h$ does not change sign when $\lambda>0$, in particular when $H\ge0$, will help us prove that for some choices of $C$ the hypersurface $M$ is embedded.

\item If we assume without loss of generality that

$$\frac{1}{|B_1(0)|}\, B_1(0)=(0,\dots,1,0,0),\, B_2(0)=(0,\dots,0,1,0)\com{and}
B_3(0)=(0,\dots,0,0,1)$$

then,

\begin{eqnarray*}
B_1(u)&=&\frac{1}{r}(0,\dots0,1,0,0)\\
B_2(u)&=&\sin(\theta(u))(0,\dots0,0,1)+\cos(\theta(u)) (0,\dots,0,1,0)\\
B_3(u)&=&\cos(\theta(u))(0,\dots0,0,1)-\sin(\theta(u)) (0,\dots,0,1,0)
\end{eqnarray*}

where $\theta:{\bf R}\to {\bf R}$ is a given by

$$\theta(u)=\int_0^u\frac{r(s)\lambda(s)}{1-r^2(s)}ds$$
\item If

$$K=K(H,n,C)=\theta(T)=\int_0^T\frac{r(s)\lambda(s)}{1-r^2(s)}ds=2\, \int_0^{\frac{T}{2}} \frac{r(s)\lambda(s)}{1-r^2(s)}ds$$

then, for any positive integer $m$ and any $u\in [mT,(m+1)T]$ we have that

$$\theta(u)=mK+\theta(u-mT)$$

This property is a consequence of the existence and uniqueness theorem for differential equation and will be used to prove the invariance of $M$ under some rotations.

 \item If ${q(u)=\gamma(u)+r^2(u)\eta(\gamma(u))}$, then
 $$\<q,q\>=1-r^2\com{and }B_2=\frac{q}{|q|}\com{i.e} q=\sqrt{1-r^2}\, B_2$$
\end{enumerate}

\subsection{A classification of constant mean curvature hypersurfaces in spheres with two principal curvatures }

We are ready to define the examples of constant mean curvature hypersurfaces on $S^{n+1}$ when $n\ge2$. Here is the theorem:

\begin{thm}\label{construction Thm}
Let $n$ be a positive integer greater than $1$ and let $H$ be a non negative real number.

\begin{enumerate}

\item
Let $g_C:{\bf R}\to {\bf R}$ be a $T$-periodic solution of the equation (\ref{integral equation for w 1}) associated with this $H$ and a positive constant $C$. If $\lambda,r,\theta:{\bf R}\to {\bf R}$ are defined by

$$r=\frac{g_C}{\sqrt{C}},\quad  \lambda=H+g_C^{-n} \com{and}\theta(u)=\int_0^u\frac{r(s)\lambda(s)}{1-r^2(s)}ds$$

then, the map $\phi:S^{n-1}\times {\bf R}\to S^{n+1}$ given by

\begin{eqnarray}\label{the immersions}
\phi(y,u)=(\, r(u)\, y,\sqrt{1-r(u)^2}\, \cos(\theta(u)),\sqrt{1-r(u)^2} \, \sin(\theta(u))\, )
\end{eqnarray}

is an immersion with constant mean curvature $H$.

\item

If $K(H,n,C)=2\int_0^\frac{T}{2}\frac{r(u)\lambda(u)}{1-r^2(u)}du=\frac{2 \pi}{k}$ for some positive integer $k$, then, the image of the immersion $\phi$ is an embedded compact hypersurface in $S^{n+1}$. In general, we have that if $K(H,n,C)=\frac{2 k \pi}{m}$ for a pair of integers, then, the image of the immersion $\phi$ is a compact hypersurface in $S^{n+1}$.

\item

 Let $n$ be an integer greater than $2$, and let $M\subset S^{n+1}$ be a connected compact hypersurface with two principal curvatures $\lambda$ with multiplicity $n-1$, and $\mu$ with multiplicity 1. If $\lambda-\mu$ is positive and the mean curvature $H=(n-1)\lambda+\mu$ is non negative and constant, then, up to a rigid motion of the sphere, $M$ can be written as an immersion of the form (\ref{the immersions}). Moreover, $M$ contains in its group of isometries the group $O(n)\times Z_m$, $m$ is the positive integer such that $K(H,n,C)=\frac{2 k \pi}{m}$, with $k$ and $m$ relative primes.

\end{enumerate}

\end{thm}

\begin{proof}

Defining $B_1$ and $B_2$ as before we have that

$$ \phi(y,u)=r(u)(y,0,0)+\sqrt{1-r(u)^2}B_2(u) $$

A direct verification shows that,

$$\frac{\partial \phi}{\partial u}=r^\prime\, (y,0,0)-\frac{r\, r^\prime}{\sqrt{1-r^2}}\, B_2
+\frac{\lambda \, r}{\sqrt{1-r^2}}\, B_3$$

We have that $\<\frac{\partial \phi}{\partial u} ,\frac{\partial \phi}{\partial u}\>=1$ and that the tangent space of the immersion at $(y,u)$ is given by

$$T_{\phi(y,u)}=\{(v,0,0)+s\, \frac{\partial \phi}{\partial u}: \<v,y\>=0\com{and} s\in{\bf R}\}$$

A direct verification shows that the map

$$\nu(y,u)=-r(u)\lambda(u)\,(y,0,0)+\frac{r^2(u)\, \lambda(u)}{\sqrt{1-r^2(u)}}\, B_2(u)+\frac{r^\prime(u)}{\sqrt{1-r^2(u)}}\, B_3(u)$$

satisfies that $\<\nu,\nu\>=1$, $\<\nu,\frac{\partial \phi}{\partial u}\>=0$ and for any $v\in\bfR{n}$ with $\<v,y\>=0$ we have that $\<\nu,(v,0,0)\>=0$. It then follows that $\nu$ is a Gauss map of the immersion $\phi$. The fact that the immersion $\phi$ has constant mean curvature $H$ follows because for any unit vector $v$ in $\bfR{n}$ perpendicular to $y$, we have that

$$\beta(t)=(r\cos(t)y+r\sin(t)v,0,0)+\sqrt{1-r^2}\, B_2=\phi(\cos(t)y+r\sin(t)v,u)$$

satisfies that $\beta(0)=\phi(y,u)$, $\beta^\prime(0)=rv$ and

$$\frac{d \nu(\beta(t))}{dt}\big{|}_{t=0} = d\nu(rv)=-r\lambda\, v$$

Therefore, the tangent vectors of the form $(v,0,0)$ are principal directions with principal curvature $\lambda$ and multiplicity $n-1$. Now, since $\<\frac{\partial \phi}{\partial u},(v,0,0)\>=0$, we have that $\frac{\partial \phi}{\partial u}$ defines a principal direction, i.e. we must have that $\frac{\partial \nu}{\partial u}$ is a multiple of $\frac{\partial \phi}{\partial u}$.  A direct verification shows that if we define $\mu:{\bf R}\to {\bf R}$ by $\mu(u)=nH-(n-1)\lambda(u)$, then,

$$\<\frac{\partial \nu}{\partial u},y\>=-\lambda^\prime \, r-\lambda r^\prime=(\lambda-\mu)r^\prime-\lambda\, r^\prime=-\mu\, r^\prime=-(nH-(n-1)\lambda)r^\prime$$

We also have that $\<\frac{\partial \phi}{\partial u},y\>=r^\prime$, therefore,

$$\frac{\partial \nu}{\partial u}= d\nu(\frac{\partial \phi}{\partial u})=-\mu\, \frac{\partial \phi}{\partial u}=-(nH-(n-1)\lambda) \frac{\partial \phi}{\partial u}$$

It follows that the other principal curvature is $nH-(n-1)\lambda$. Therefore $\phi$ defines an immersion with constant mean curvature $H$, this proves the first item in the Theorem.

In order to prove the second item,  we notice that if $K(H,n,C)=\frac{2\pi}{k}$ for some positive $k$, then, we get that $\theta(kT)=2\pi $, this last fact  makes the image of the immersion $\phi$ compact. This immersion is  embedded because the immersion $\phi$ is one to one for values of $u$ between $0$ and $k T$ as we can easily check using the fact that  whenever $H\ge0$, the function $\theta$ is strictly increasing. Recall that under these circumstances
$\theta(0)=0$ and $\theta(kT)=2\pi$. The prove of the other statement in this item is similar.

Let us prove the next item. For $n>2$, let us now consider a minimal hypersurface $M$ with the properties of the statement. We will use the notation that we used in the preliminaries, in particular the function $w:M\to {\bf R}$ is defined by the relation $(\lambda-\mu)=nw^n$, in particular we will assume that $B_1(0)$, $B_2(0)$ and $B_3(0)$ are chosen as before. By Theorem (\ref{equations}) we get that the distribution
$\hbox{Span}\{e_1,\dots,e_{n-1}\}$ is completely integrable. Let us fix a point $p_0$ in $M$ and let us define the geodesic $\gamma:{\bf R}\to M$, and the functions $r:{\bf R}\to {\bf R}$ as before and let us denote by $M_u\subset M$ the $(n-1)$-dimensional integral submanifold of $M$ of this distribution that passes through $\gamma(u)$. Let us define the vector field $\eta$ on $M$ as before. Recall that  $B_1(u)=\eta(\gamma(u))$. Fixing a value $u$, let us define the maps

$$\rho_u,\zeta_u:M_u\to \bfR{{n+2}}\com{by} \rho_u(x)=x+\frac{w^2(x)}{C}\,\eta\com{and} \zeta_u(x)=\nu(x)+\lambda(x)\, x$$

Using the equations in section (\ref{properties of eta}) we get that the maps $\rho_u$ and $\zeta_u$ are constant. Therefore,

$$\rho_u(x)=x+\frac{w^2(x)}{C}\,\eta=\gamma(u)+r^2(u)B_1=\sqrt{1-r^2} B_2$$

Notice that for every $x\in M_u$, we have that

$$|x-\rho_u(x)|^2=|Z(u)-\sqrt{1-r^2} B_2(u)|^2=r^2(u)$$

Therefore $M_u$ is contained in a sphere with center in $\sqrt{1-r^2} B_2$ and radius $r$.
We have that the vectors $e_1,\dots,e_{n-1}$ are perpendicular to the vectors

$$\rho_u(x)=\sqrt{1-r^2(u)} B_2(u)\com{and} \zeta_u(x)=Y(u)+\lambda(u) Z(u)$$

Since $\<Y(u)+\lambda(u) Z(u),B_1(u)\>=0$, $\<Y(u)+\lambda(u) Z(u),B_2(u)\>=\frac{\lambda r^2}{\sqrt{1-r^2}}$ and $\<Y(u)+\lambda(u) Z(u),B_3(u)\>=\frac{r^\prime}{\sqrt{1-r^2}}$, we get that

$$\zeta_u(x)=\frac{\lambda r^2}{\sqrt{1-r^2}} \, B_2+ \frac{r^\prime}{\sqrt{1-r^2}}\, B_3$$

It follows that, anytime $r^\prime(u)\ne 0$, we have that  all tangent vectors of $M_u$ lies in the $n$-dimensional space perpendicular to the two dimensional space spanned by $B_1(u)$ and $B_2(u)$.  Since this two dimensional space is independent of $u$, we conclude that every point  $x\in M_u$, satisfies that

$$ x-\rho_u(x)=r(u)(y,0,0)\com{where} |y|^2=1$$

or equivalently,

$$x=r(u)\, (y,0,0)+\rho_u(x)=r(u)(y,0,0)+\sqrt{1-r(u)^2}B_2(u)$$

Since the set of points where $r^\prime$ is discrete, we conclude that the expression for the points $x\in M_u$ holds true for all $u$. The theorem then follows because the manifold $M$ is connected.

The property on the group of isometries of the manifold follows because we can write  $M$ as the image of the map

\begin{eqnarray}
\phi(y,u)=(\, r(u)\, y,\sqrt{1-r(u)^2}\, \cos(\theta(u)),\sqrt{1-r(u)^2} \, \sin(\theta(u))\, )
\end{eqnarray}

The group $O(n)$ is part of the isometries of $M$ because any isometry in $\bfR{{n+2}}$ that fixes the origin and the last two entries of $\bfR{{n+2}}$ leave the our manifold $M$ invariant. The group $Z_m$ is included in the groups of isometries because the close curve given by the last two entries is built by gluing $m$ pieces of the the curve

$$\alpha(u)=(\sqrt{1-r(u)^2}\, \cos(\theta(u)),\sqrt{1-r(u)^2} \, \sin(\theta(u)))\quad 0\le u\le K(H,n,C)=\frac{2 k \pi}{m}$$

This last statement is true by the the following observation already pointed out in the previous section.

$$ \com{For any positive integer $j$ and $u\in [jT,(j+1)T]$ we have that }\theta(u)=jK+\theta(u-jT) $$

\end{proof}

\begin{cor}
If $M$ is one of the examples in the previous theorem with $H=0$, then, the stability index, i.e, the number of negative eigenvalues of the operator $J(f)=-\Delta f-nf-||A||^2\, f$ is greater than $2n+5$.

\end{cor}

\begin{proof}
This follows Theorem (3.1.1) in (\cite{P1}) and the fact that the set of isometries of  $M$ contains a subgroup that moves every point in $M$.
\end{proof}

\section{Embedded CMC surfaces in $S^3$}

For surfaces in $S^3$, the examples will have the following form

$$\Sigma= \{(\, \sqrt{1-|\alpha(t)|^2}\, \cos(s),\sqrt{1-|\alpha(t)|^2}\, \sin(s),\alpha(t)\,) \, :\,  t\in{\bf R}\com{and} s\in{\bf R}\, \}$$

where $\alpha(t)=(x_3(t),x_4(t))$ will be called the profile curve and we will refer to the surface $\Sigma$ as the rotation of the curve $\alpha$.

An interesting observation is that if we rotate the set of points in a line that lies in the unit disk, we get a totally umbilical sphere. We can see this by noticing that this surface can be written as

 $$\{(x_1,x_2,x_3,x_4): x_1^2+x_2^2+x_3^2+x_4^2=1\com{and} a x_3+b x_4=c\, \}\com{where $a$, $b$ and $c$ are constants}$$


 From this observation it follows that if we rotate a regular polygonal with $m$ sides and vertexes in the unit circle, then we obtain $m$ kissing spheres in $S^3$.   The family of constructed examples  can be subdivide into families $F_0$,
$F_2$, $F_3$ and so on, where $F_0$ are the isoparametric examples and $F_m$ are non isoparametric examples which group of isometries contains the subgroup of isometries $O(2)\times Z_m$. The family $F_m$ is a collections of surfaces that moves from the rotation of an regular polygon with $m$ sides to the rotation of a circle. In this way, the surfaces in this paper resemble Delaunay's CMC surfaces in $\bfR{3}$, in the sense that the examples here, move from a CMC  without singularities (the rotation of a circle, i.e. an isoparametric CMC surface), to a CMC with $m$ singularities (the union of $m$ spheres). Recall that Delaunay's surfaces move from a constant mean curvature without singularities (a cylinder, the rotation of a line), to a constant mean curvature with infinitely many singularities (a collection of round spheres, the rotations of a collection of semi-circles).

In the case $n=2$ we can explicitly compute the roots, $t_1$ and $t_2$, of the polynomial given in (\ref{the polynomial}), therefore, if $H$ is any real number, and $C$ is a constant greater than $2(H+\sqrt{1+H^2})$, then the solution $g:{\bf R}\to {\bf R}$ that satisfies the equation (\ref{integral equation for w 1}), i.e. the equation

\begin{eqnarray}\label{integral equation of g n=2}
 (g^\prime)^2+g^{-2}+(1+H^2)g^2+2H=C
 \end{eqnarray}

is positive and more precisely, it varies from $t_1$ to $t_2$, where

\begin{eqnarray}\label{zeroes of the polynomial}
t_1=\sqrt{\frac{C-2H-\sqrt{C^2-4 H C-4}}{2(1+H^2)}}\com{and}
t_2=\sqrt{\frac{C-2H+\sqrt{C^2-4 H C-4}}{2(1+H^2)}}
\end{eqnarray}

Also, we have that $g$ is a periodic function with period

$$T=2 \int_{t_1}^{t_2}\frac{t}{\sqrt{(C-2H) t^2-1-(1+H^2)t^4}}dt$$

the function $\theta:{\bf R}\to {\bf R}$ previously defined is given by

$$\theta(u)=\int_0^u\frac{r(s)\lambda(s)}{1-r^2(s)}ds$$

where $r(u)=\frac{g(u)}{\sqrt{C}}$ and $\lambda(u)=1+g^{-2}(u)$.
Moreover, Theorem (\ref{construction Thm}) give us that the following map,

$$\phi(v,u)=(r(u)\cos(v),r(u)\sin(v),
\sqrt{1-r^2(u)}\cos(\theta(u)),\sqrt{1-r^2(u)}\sin(\theta(u)))$$

defines a surface with constant mean curvature $H$.

In this section we will study our constructions for the case $n=2$ and by studying the dependence of the function $r$ in term of the constant $C$ we will prove the existence of compact embedded CMC surfaces in $S^3$ with any prescribed value $H$ different from $0,\frac{1}{\sqrt{3}}$ and $-\frac{1}{\sqrt{3}}$. As pointed out in item (2) of Theorem (\ref{construction Thm}), in order to guarantee the existence of embedded surfaces we need to study the number $K(H,2,C)$ where,

\begin{eqnarray}\label{K for n equal 2}
K=K(H,2,C)=2\int_0^{\frac{T}{2}}\frac{r(s)\lambda(s)}{1-r^2(s)}ds
\end{eqnarray}

The following lemma help us to understand $K$.

\begin{lem}\label{a nice expression for K}
For any $C>2\sqrt{1+H^2}+2H$, the number $K(H,2,C)$ given in (\ref{K for n equal 2}) can be computed with the following integral

$$K(H,2,C)=\int_0^{\pi}\frac{\frac{1}{C}+H(-q_2\cos(t)-q_1)}{(1+q_1+q_2\cos(t))\sqrt{-q_2\cos(t)-q_1}\sqrt{1+H^2}}dt$$

where

$$q_1=\frac{-C+2H}{2 C+2 C H^2}\com{and} q_2=\frac{\sqrt{-4+C^2-4 C H}}{2 (C+ H^2) }$$

Moreover for any $H> 0$,

$$K(H,2,C)\to \frac{\sqrt{2} \, \pi\, (H+\sqrt{1+H^2})^\frac{3}{2}}{(1+H^2)^{\frac{1}{4}}(1+2H^2+2 H\sqrt{1+H^2})}\com{when} C\to 2(H+\sqrt{1+H^2})^+$$
and
$$K(H,n,C)\to 2 \hbox{\rm ArcCot}(H)\com{when} C\to \infty$$

\end{lem}

\begin{proof}
Notice that the function $r$ is strictly increasing in $[0,\frac{T}{2}]$ and that $r(0)=\frac{t_1}{\sqrt{C}}$ and $r(\frac{T}{2})=\frac{t_2}{\sqrt{C}}$. After we make the change of variable $t=r(s)$, using (\ref{equation for r}), the expression for $K(H,2,C)$ in  (\ref{K for n equal 2}) reduces to

$$K(H,2,C)=2\int_{\frac{t_1}{\sqrt{C}}}^{\frac{t_2}{\sqrt{C}}} \frac{t(H+c^{-1} t^{-2})}{(1-t^2)\sqrt{1-t^2(1+H^2+2HC^{-1}\, t^{-2}+C^{-2}t^{-4})}}\, dt$$

Now, making $u=t^2$ we get that

\begin{eqnarray*}
K(H,2,C)&=&\int_{\frac{t_1^2}{{C}}}^{\frac{t_2^2}{{C}}} \frac{(H+c^{-1} u^{-1})}{(1-u)\sqrt{1-u(1+H^2+2HC^{-1}\, u^{-1}+C^{-2}u^{-2})}}\, du\\
&=&\int_{\frac{t_1^2}{{C}}}^{\frac{t_2^2}{{C}}} \frac{(H u+c^{-1})}{\sqrt{u}(1-u)\sqrt{-u^2(1+H^2)+(1-2HC^{-1})\, u - C^{-2})}}\, du
\end{eqnarray*}

The expression inside the radical in the integral above can be written as

$$-u^2(1+H^2)+(1-2HC^{-1})\, u - C^{-2}=(1+H^2)\big{(}\,q_2^2-(u+q_1)^2\,\big{)}$$

where,

$$q_1=\frac{2H-C}{(1+H^2) 2C}\com{and} q_2=\frac{\sqrt{C^2-4 C H-4}}{2 C (1+H^2)} $$

By doing the substitution $v=u+q_1$ we obtain that
\begin{eqnarray*}
K(H,2,C)=\int_a^b \frac{(H (v-q_1)+c^{-1})}{\sqrt{1+H^2}\, \sqrt{v-q_1}\, (1-v+q_1)\, \sqrt{q_2^2-v^2}}\, dv
\end{eqnarray*}

where

$$b=\frac{t_2^2}{C}+q_1=q_2\com{and}a=\frac{t_1^2}{C}+q_1=-q_2$$
The integral formula in the lemma follows by doing the substitution $v=-q_2\cos(t)$.
The formula for the  limit of $K(H,2,C)$ follows by taking the limit of the function inside the integral and then integrating.

\end{proof}

{\bf Remark}

{\sl
For the case $H=0$, Otsuki in (\cite{O}), showed that $K\to \pi$ when $C\to \infty$ and $K\to \sqrt{2}\, \pi$ when $C\to 2$. These values agree with the limit when $H\to 0$ of the expressions that we have.}

\begin{thm}

For any $H$ different from $0$ and $\pm\frac{1}{\sqrt{3}}$ there exists a non isoparametric  embedded compact surface with constant mean curvature $H$ in $S^3$.

\end{thm}

\begin{proof}
Without loss of generality we will assume that $H>0$. By Theorem (\ref{construction Thm}) we only need to check that for every $H\ne \frac{1}{\sqrt{3}}$, there exist a positive integer $m$ and a constant $C$ such that the number $K(H,2,C)=\frac{2\pi}{m}$. By Lemma (\ref{a nice expression for K}) we only need to prove that for some integer $m$, we have that

$$b_1(H)= 2 \hbox{\rm ArcCot}(H)<\frac{2\pi}{m}< \frac{\sqrt{2} \, \pi\, (H+\sqrt{1+H^2})^\frac{3}{2}}{(1+H^2)^{\frac{1}{4}}(1+2H^2+2 H\sqrt{1+H^2})}=b_2(H)$$

A direct computation shows that,
$$   b_2^\prime(H)=-\frac{\pi\, \sqrt{H+\sqrt{1+H^2}}}{\sqrt{2}(1+H^2)^\frac{5}{4}}$$

therefore both of the bound functions  $b_1(H)$ and $b_2(H)$ are decreasing, positive and they have $0$ as a horizontal asymptote. The function $b_1(H)= 2 \hbox{\rm ArcCot}(H)$ starts with the value $\pi$ at $H=0$ and the function $b_2(H)$ starts at $\sqrt{2}\pi$ when $H=0$. Therefore, for values of $H$ close to $0$, as long as $b_2(H)>\pi$, we can pick $m=2$ and get an embedded CMC surface with the given value $H$. Since the solution of the equation

 $$ \frac{\sqrt{2} \, \pi\, (H+\sqrt{1+H^2})^\frac{3}{2}}{(1+H^2)^{\frac{1}{4}}(1+2H^2+2 H\sqrt{1+H^2})}=b_2(H)=\pi=\frac{2\pi}{m}\com{with} m=2$$

 is  $H=\frac{1}{\sqrt{3}}$, we get that for every $H$ between $0$ and $\frac{1}{\sqrt{3}}$ there exist an embedded $H$-CMC surface.  But, is there a surface with constant mean curvature  $H=\frac{1}{\sqrt{3}}$ associated with $m=3$?. We could obtain  this if we
  have that $b_1(\frac{1}{\sqrt{3}})<\frac{2\pi}{3}$, but curiously,

 $$b_1(\frac{1}{\sqrt{3}})=2 \hbox{\rm ArcCot}(\frac{1}{\sqrt{3}})=\frac{2\pi}{3}$$

 Therefore, we can not guarantee the existence of a embedded surface with constant mean curvature $\frac{1}{\sqrt{3}}$. For values of $H$ after $\frac{1}{\sqrt{3}}$ we have the existence of embedded surfaces with constant mean curvature $H$ associated with $m=3$, as long as $b_2(H)>\frac{2 \pi}{4}$. A direct computation shows that $b_2(\frac{7}{4\sqrt{2}})=\frac{2 \pi}{4}$, therefore all values of $H$ between $ \frac{1}{\sqrt{3}}$ and  $\frac{7}{4\sqrt{2}}$ can be realized as constant mean curvatures of surfaces associated with $m=3$, but, is there a surface with constant mean curvature  $H=\frac{7}{4\sqrt{2}}$ associated with $m=4$? We could obtain  this if we
  have that $b_1(\frac{7}{4\sqrt{2}})<\frac{2\pi}{4}$. In this case we have that indeed
  $b_1(\frac{7}{4\sqrt{2}})<\frac{2\pi}{4}$ and therefore $H=\frac{7}{4\sqrt{2}}$ is the mean curvature of a CMC surfaces associated with $m=4$. We will prove the existence of embedded surfaces with constant mean curvature greater than $\frac{1}{\sqrt{3}}$ by showing that

\begin{eqnarray}\label{the inequality on the periods n=2}
b_2(\cot(\frac{\pi}{(m+1)}))>\frac{2\pi}{m}\com{for any integer} m\ge3
\end{eqnarray}
  The reason the statement above suffices is because, for any $H>\frac{1}{\sqrt{3}}$ there exists an integer $m\ge 3$ such that

  $$\cot(\frac{\pi}{m})<H\le \cot(\frac{\pi}{m+1})$$

Since $b_1$ and $b_2$ are decreasing, we get that

$$\frac{2 \pi}{m+1} \le b_1(H)<\frac{2\pi}{m}\com{and} b_2(\cot(\frac{\pi}{m+1}))<b_2(H) $$

 Using the inequality (\ref{the inequality on the periods n=2}) we get that

 $$ b_1(H)<\frac{2\pi}{m}<b_2(\cot(\frac{\pi}{m+1}))<b_2(H)$$

 which guarantees that there exists a surface with constant mean curvature $H$  that contains the group $O(2)\times Z_m$ in its group of isometries.

To prove the inequality (\ref{the inequality on the periods n=2}) it is enough to see that

$$d(H)=b_2(\cot(\frac{\pi}{(m+1)}))- \frac{2\pi}{m}=\frac{\pi \cot^\frac{3}{2}(\frac{\pi}{2m+2})\, \sec^4(\frac{\pi}{2m+2})}{2\sqrt{2}\csc^\frac{5}{2}(\frac{\pi}{m+1})}-\frac{2\pi}{m}$$

and that

$$d^\prime(H)=-\pi \frac{m^2 \, \pi \cos(\frac{\pi}{2m+2})-2(1+m)^2}{m^2(1+m)^2}$$

For $m\ge 3$ is not difficult to see that $m^2 \, \pi \cos(\frac{\pi}{2m+2})-2(1+m)^2$ is positive. Since $d(3)>0$, the limit when $m\to \infty$ is zero and $d$ is strictly decreasing, then we conclude that $d$ must be positive. Therefore, the theorem follows.
\end{proof}

\begin{cor}
For any $H\in (0,\frac{1}{\sqrt{3}})$, there exists a $O(2)\times Z_2$-invariant surface with constant mean curvature $H$. For any $H\in (\frac{1}{\sqrt{3}},\frac{7}{4\sqrt{2}})$ there exists a $O(2)\times Z_3$-invariant surface with constant mean curvature $H$. For any $H\in (1,\frac{7}{\sqrt{15}})$ there exists a $O(2)\times Z_4$-invariant surface with constant mean curvature $H$. In general for any $m>1$, for any $H\in (b_1^{-1}(\frac{2 \pi}{m}),b_2^{-1}(\frac{2 \pi}{m}))$ there exists a $O(2)\times Z_m$-invariant surface with constant mean curvature $H$. Notice that for values of $H$ between $1$ and $\frac{7}{4\sqrt{2}}$ there are surface with constant mean curvature $H$ invariant under  $O(2)\times Z_3$ and $O(2)\times Z_4$. The same overlapping occurs for any $m>2$.

\end{cor}

\section{Hypersurface with CMC in $S^{n+1}$, general case}

In this section we will study the existence of compact examples in $S^{n+1}$ by studying the values $K(H,n,C)$. The key lemma in this study is the following.

\begin{lem}\label{the lemma}
Let $f:(-\delta,\delta)\to {\bf R}$ be a smooth function such that $f(0)=f^\prime(0)=0$ and $f^{\prime\prime}(0)=-2a<0$. For positive values of $c$ close to $0$, let $t(c)$ be the first positive  root of the function $f(t)+c$. We have that

$$\int_0^{t(c)}\frac{dt}{\sqrt{f(t)+c}}=\frac{\pi}{2 \sqrt{a}}$$
\end{lem}

\begin{proof}

For any $b>a$ let us define the function $h(t)=f^\prime(t)+2bt$. Since $h^{\prime}(0)=2(b-a)>0$ there exists a positive $\epsilon$ such that $h^\prime(t)>0$ for all $t\in [0,\epsilon]$. Now for any $c$ such that $t(c)<\epsilon$ the function

$$g(t)=f(t)+c-(bt(c)^2-bt^2)$$

satisfies that $g(t(c))=0$ and $g^\prime(t)=h(t)>0$. Therefore, $g(t)<0$ for any $t\in[0,t(c)]$. By the definition of $g(t)$ we get that

$$0<f(t)+c<bt(c)^2-bt^2\com{for all } t\in[0,t(c))$$

and therefore we get that

$$\frac{\pi}{2\sqrt{b}}=\int_0^{t(c)}\frac{dt}{\sqrt{bt(c)^2-bt^2}}<
\int_0^{t(c)}\frac{dt}{\sqrt{f(t)+c}}$$

Likewise, for any $b<a$, the same argument shows that

$$\int_0^{t(c)}\frac{dt}{\sqrt{f(t)+c}}<
\int_0^{t(c)}\frac{dt}{\sqrt{bt(c)^2-bt^2}}=\frac{\pi}{2\sqrt{b}}$$

Since $b\ne a$ can be chosen arbitrarily close to $a$, we conclude the lemma.

\end{proof}

\begin{cor}\label{cor of the lemma}
Let $\epsilon$ and $\delta$ be positive real numbers and
let $f:(t_0-\epsilon,t_0+\epsilon)\to {\bf R} $ and $g:(-\delta,\delta)\times (t_0-\epsilon,t_0+\epsilon)\to {\bf R}$ be two smooth functions such that $f(t_0)=f^\prime(t_0)=0$ and $f^{\prime\prime}(t_0)=-2a<0$. If for any small $c>0$, $t_1(c)<t_0<t_2(c)$ are such that $f(t_1(c))+c=0=f(t_2(c)))+c$, then

$$\lim_{c\to0^+}\int_{t_1(c)}^{t_2(c)}\frac{g(c,t)\, dt}{\sqrt{f(t)+c}}=\frac{g(0,t_0)\, \pi }{ \sqrt{a}}$$

\end{cor}

This lemma allows us to prove the main theorem in this paper,

\begin{thm}
For any $n\ge 2$  and any  $H\in(0,\frac{2\sqrt{n-1}}{n\sqrt{3}})$ there exists a non isoparametric compact embedded hypersurface in $S^{n+1}$ with constant mean curvature $H$. More generally, for any integer $m>1$ and $H$ between the numbers

$$\cot{\frac{\pi}{m}} \com{and} \frac{ (m^2-2)\, \sqrt{(n-1)}}{n \sqrt{m^2-1}} $$

 there exist a non isoparametric compact embedded hypersurface in $S^{n+1}$ with constant mean curvature $H$ such that its group of isometries contains the group $O(n)\times Z_m$.

\end{thm}

\begin{proof}

We will consider only positive values for $H$. Here we will use the explicit solution for the ODE (\ref{integral equation for w 1}) given in section (\ref{section the function w and its solutions}). Let us rewrite this ODE as,

$$ (g^\prime)^2=q(g) \com{where} q(v)=C-v^{2-2n}-(1+H^2)v^2-2Hv^{2-n}$$

We already pointed out in section (\ref{section the function w and its solutions}) that for some values of $C$, the function $q$ has positive values between two positive roots of $q$, denoted by $t_1$ and $t_2$. Let us be more precise and give an expression for how big $C$ needs to be. A direct verification shows that

$$q^\prime(v)=-2 (1 + H^2) v - (2 - 2 n) v^{1 - 2 n} - 2 H (2 - n) v^{1 - n}$$

  and that the only positive root of $q^\prime$ is

 \begin{eqnarray}\label{def v0}
 v_0=(\frac{\sqrt{H^2n^2+4(n-1)}+(n-2) H}{2+2 H^2})^{\frac{1}{n}}
  \end{eqnarray}
  Therefore, for positive values of $v$, the function $q$ increases from $0$ to $v_0$ and
  decreases for values greater than $v_0$.  A direct computation shows that $q(v_0)= c-c_0$
  where,

 \begin{eqnarray}\label{def c0}
 c_0=n\, (2+2 H^2)^{\frac{n-2}{n}}\,
 \frac{2+n H^2 + H\sqrt{H^2n^2 + 4 (n-1)}}{\big{(}\, (n-2)H +\sqrt{H^2n^2+ 4 (n-1)}\, \big{)}^{\frac{2n-2}{n}} }
   \end{eqnarray}

 Therefore, whenever $C>c_0$ we will have exactly two positive roots of the function $q(v)$ that we will denote by $t_1(C)$ and $t_2(C)$ to emphasize its dependents on $C$. A direct computation shows that $q^{\prime\prime}(v_0)= -2 a$ where

 $$a=2n(1+H^2)\, \frac{4(n-1)+H^2\, n^2+H\, (n-2)\, \sqrt{4(n-1)+H^2 n^2}}{(\, H(n-2)+\sqrt{4(n-1)+H^2 n^2}\, )^2}$$

 Using the notation and results of section (\ref{section vector fields that lie on a planes}), we get that
\begin{eqnarray}\label{expression for K(H,n)}
K(H,n,C)=2\int_0^{\frac{T}{2}}\frac{r(s)\lambda(s)}{1-r^2(s)}\, ds
\end{eqnarray}
Since $r(s)=\frac{g(s)}{\sqrt{C}}$ and $\lambda(s)=H+g(s)^{-n}$ we get that

$$K(H,n,C)=2\int_0^{\frac{T}{2}}\frac{\sqrt{C}g(s)(H+g^{-n}(s))}{c-g^2(s)}\, ds$$

Since $g(0)=t_1(c)$ and $g(\frac{T}{2})=t_2(c)$, by doing the substitutions $t=g(s)$ we get

$$K(H,n,C)=2\int_{t_1(c)}^{t_2(c)}\frac{ \sqrt{C}t(H+t^{-n})}{c-t^2} \frac{1}{\sqrt{q(t)}}  \, dt  $$

 Since $a>0$ we can apply Corollary (\ref{cor of the lemma}) to the get that

 $$\lim_{c\to c_0^+}K(H,n,C)=\pi\, \sqrt{2-\frac{2 n H}{\sqrt{4(n-1)+H^2n^2}}}$$

It can be verified that  this bound is the same bound we found for the case $n=2$.

In order to analyze the limit of the function $K(H,n,C)$ when $C\to \infty$ we return to the expression (\ref{expression for K(H,n)}) and we make the substitution  $t=r(s)$ to obtain

$$K(H,n,C)=2\int_{\frac{ t_1(C)}{\sqrt{C}}}^{\frac{ t_2(C)}{\sqrt{C}}} \frac{t (H+C^{-\frac{n}{2}}\, t^{-n})}{(1-t^2)\, \sqrt{1-t^2(1 +(H+C^{-\frac{n}{2}}\, t^{-n})^2)}}\, dt$$

In this case we have used the equation (\ref{equation for r}) to change the $ds$ for the $dt$. Notice that the limit values $\frac{ t_1(C)}{\sqrt{C}}$ and $\frac{ t_2(C)}{\sqrt{C}}$ can also be characterize as the only positive roots of the function,

$$ \tilde{q}=1-t^2(1 +(H+C^{-\frac{n}{2}}\, t^{-n})^2)=1-(1+H^2)t^2-C^{-n}t^{2-2n}-2 H C^{-\frac{n}{2}}t^{2-n}$$

because of the relation $q(v)=C \tilde{q}(\frac{v}{\sqrt{C}})$. Since $\tilde{q}(\frac{1}{\sqrt{1+H^2}})<0$ and for every positive $\epsilon<\frac{1}{\sqrt{1+H^2}}$ we have that

$$\lim_{C\to\infty}\tilde{q}(\epsilon)>0 \com{and} \lim_{c\to\infty}\tilde{q}(\frac{1}{\sqrt{1+H^2}}-\epsilon)>0$$

then, we conclude that the only two positive roots of $\tilde{q}$ converge to $0$ and to $\frac{1}{\sqrt{1+H^2}}$ when $C\to \infty$. Therefore,

$$\lim_{C\to\infty}K(H,n,C)=2 \int_0^{\frac{1}{\sqrt{1+H^2}}}\frac{H t}{(1-t^2)\, \sqrt{1-(1+H^2)t^2}}\, dt= 2 \hbox{\rm ArcCot}(H) $$

Notice that this bound is the same bound we found for the case $n=2$. Therefore, for any fixed $H>0$, the function $K(H,n,C)$ takes all the values between

$$a_1(H)=2\hbox{\rm ArcCot}(H)\com{and} a_{2,n}(H)=\pi\, \sqrt{2-\frac{2 n H}{\sqrt{4(n-1)+H^2n^2}}}$$

 We have that the functions $a_1(H)$ and $a_{2,n}(H)$ are decreasing. Moreover, we have that for any $y<\sqrt{2}$

$$a_{2,n}(\frac{2\, (2-y^2)\, \sqrt{n-1}}{n\, y\, \sqrt{4-y^2}})=\pi \, y$$

Therefore, replacing $y$ by $\frac{2}{m}$ in the expression above, we obtain that for values of $H$ between

$$\cot{\frac{\pi}{m}} \com{and} \frac{ (m^2-2)\, \sqrt{(n-1)}}{n \sqrt{m^2-1}} $$

the number $\frac{2 \pi}{m}$ lies between $a_1(H)$ and $a_{2,n}(H)$, and therefore, for some constant $C$, we will have that $K(H,n,C)=\frac{2 \pi}{m}$. Using the same arguments we use in the case $n=2$ the theorem will follow. Notice that when $m=2$ these two bounds are
$0$ and $\frac{2\, \sqrt{n-1}}{n\, \sqrt{3}}$

\end{proof}

Let us finish  this section with a remark already pointed out by Otsuki in (\cite{O}).

\begin{lem}\label{Close to n}
For any integer $n\ge 2$ and any $\epsilon>0$ there exist compact non isoparametric minimal hypersurfaces in $S^{n+1}$ such that $n-\epsilon\le ||A||^2(p) \le n+\epsilon$ for all $p\in M$.
\end{lem}

\begin{proof}

This is a consequence of the fact that the expression for $v_0$ in (\ref{def v0})   reduces to $(n-1)^{\frac{1}{2n}} $ when $H=0$ and the fact that by picking $c$ close to $c_0$,  the roots $t_1(C)$ and $t_2(C)$ of the function $q$ are as close as $v_0$ as we want. Since the range of the function $g$ move from $t_1(C)$ to $t_2(C)$,  we can make the values of $g$ to move as close of $ (n-1)^{\frac{1}{2n}}$ as we want. When $H=0$, we have that

$$\lambda=g^{-n}\quad \mu=-(n-1)g^{-n} \com{and} ||A||^2= (n-1) g^{-2n}+(n-1)^2 g^{-2n}=n(n-1)g^{-2n}$$

Therefore, we can make $||A||^2$ as close $n$ as we want. By density of the rational number and the continuity of the function $K(H,n,C)$, we can choose $C$ so that $K(H,n,C)$ is of the form $\frac{2 k \pi}{m}$ for some pair of integers $m$ and $k$. This last condition will guarantee the compactness of the profile curve and therefore the compactness of the hypersurface.

\end{proof}

\section{Non isoparametric stable cones in $S^{n+1}$}

For any  compact minimal hypersurface $M\subset S^{n+1}$, let us define the operator $L_1$ and the number $\lambda_1$ as follows,

$$L_1(f)=-\Delta f-||A||^2 f\com{and}\lambda_1=\com{first eigenvalue of $L_1$}$$

 Moreover, let us denote  by $CM=\{tm:t\in [0,1],\, m\in M\, \}$ the cone over $M$. We will say that $CM$ is stable if every variation of $CM$, which holds $M$ fixed, increases area.

In (\cite{S}, Lemma 6.1.6) Simons proved that if $\lambda_1+ (\frac{n-1}{2})^2>0$ then $CM$ is stable. We will prove that for any $n\ge 6$, the cone over some non isoparametric  examples studied in this paper for $H=0$, i.e, the cone over some of the Otsuki's examples,  are stable. More precisely we have,

\begin{thm}
For any $n\ge 6$, there are non isoparametric compact hypersurfaces in $S^{n+1}$ such that their cone is stable.
\end{thm}

\begin{proof}

A direct verification shows that

$$ (\frac{n-1}{2})^2\ge n+\frac{1}{4}\com{for all} n\ge 6 $$

Using Lemma (\ref{Close to n}), let us consider a non isoparametric compact minimal hypersurface $M$ such that $||A||^2\le n+\frac{1}{8} $. We have that the first eigenvalue $\lambda_1$ of the operator $L_1$ is greater than $-n-\frac{1}{8}$ because

$$\lambda_1=\hbox{\rm inf } \{\frac{\int_M(-\Delta f-||A||^2\, f)f}{\int_Mf^2}:f\com{is smooth and $\int_Mf^2\ne0$ }\}$$

and we have that,

$$ \frac{\int_M(-\Delta f-||A||^2\, f)f}{\int_Mf^2}=\frac{\int_M |\nabla f|^2}{\int_Mf^2}-\frac{\int_M||A||^2\, f^2}{\int_Mf^2}\ge-(n+\frac{1}{8})$$

Therefore, we get that

$$\lambda_1+(\frac{n-1}{2})^2\ge -(n+\frac{1}{8})+n+\frac{1}{4}=\frac{1}{8}>0$$

which implies by Simons' result that the $CM$ is stable.

\end{proof}

\section{Some explicit solutions}

In this section we will pick some arbitrary values of $H$ to explicitly show the embedding, the graph of the profile curves, and the stereographic projections of some examples of surfaces with CMC in $S^3$.

A direct computation shows that the solution of the equation (\ref{integral equation of g n=2}) is given by

$$g(t)=\sqrt{\frac{(C - 2 H) +
 \sqrt{-4 + C^2 - 4 C H}\, \sin(2 \sqrt{1 + H^2}\, t\,)}{2 (1 + H^2))}\, }$$

From the expression for $g$ we get that its period $T$ is $\frac{\pi}{\sqrt{1+H^2}}$.
In this case, the condition on $C$ to get solutions of  the ODE  (\ref{integral equation of g n=2})  reduces to $C>2(H+\sqrt{1+H^2})$.

We can get surfaces associated with $m=2$ if we take $H$ between $0$ and
$\frac{1}{\sqrt{3}}\simeq 0.57735$ and we can surfaces associated with $m=3$ if we take $H$ between $\frac{1}{\sqrt{3}}$ and $\frac{7}{4\sqrt{2}}\simeq 1.23744$. Once we have picked the value for $H$ in the right range, in order to get the embedded surface, we need to solve the equation

$$K(H,2,C)=\int_0^{\frac{\pi}{\sqrt{1+H^2}}}\frac{\sqrt{C}\, g(t)(H+g(t)^{-2})}{C-g(t)^{2}}\, dt=\frac{2\pi}{m}$$

Finally, when we have the $H$ and the $C$, the profile curve is given by

$$ (\, \sqrt{1-\frac{g^2(t)}{C}}\, \cos(\theta(t)), \sqrt{1-\frac{g^2(t)}{C}}\, \sin(\theta(t))\, )\com{where} \theta(t)=\int_0^t\frac{\sqrt{C}\, g(\tau)(H+g^{-2}(\tau))}{C-g^2(\tau)}\, d\tau$$

and the embedding is given by

$$ (\, \frac{g(t)}{\sqrt{C}}\, \cos(u), \frac{g(t)}{\sqrt{C}}\, \sin(u), \sqrt{1-\frac{g^2(t)}{C}}\, \cos(\theta(t)), \sqrt{1-\frac{g^2(t)}{C}}\, \sin(\theta(t))\, )\quad 0\le u< 2\pi\quad 0\le t<  m\, \frac{\pi}{\sqrt{1+H^2}}$$

Here are some graphics,

\begin{figure}[h]\label{profile h=0.1}
\centerline{\includegraphics[width=10cm,height=5.5cm]{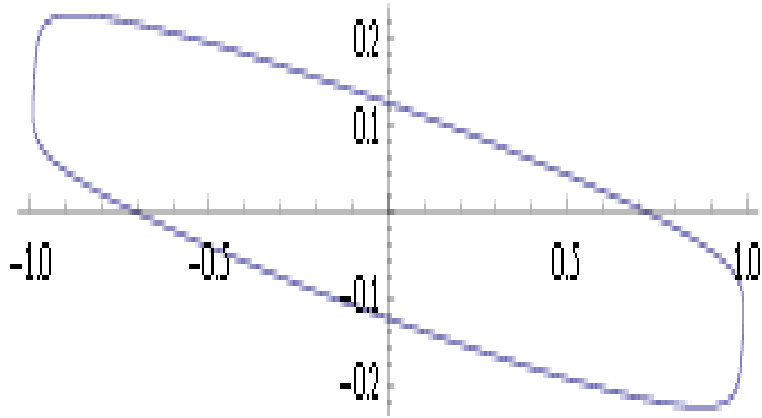}}
\caption{Profile curve for $m=2$, $H=0.1$, in this case $C=41.28796038772471$}
\end{figure}

\vfil
\vfill
\eject

\begin{figure}[h]\label{profile h=0.3}
\centerline{\includegraphics[width=10cm,height=5.5cm]{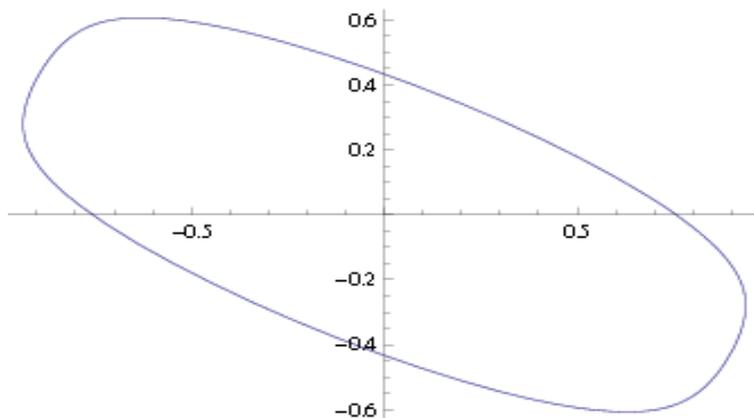}}
\caption{Profile curve for $m=2$, $H=0.3$, in this case $C=9.129645968138256$}
\end{figure}

\vskip2cm

\begin{figure}[h]\label{profile h=0.57}
\centerline{\includegraphics[width=10cm,height=5.5cm]{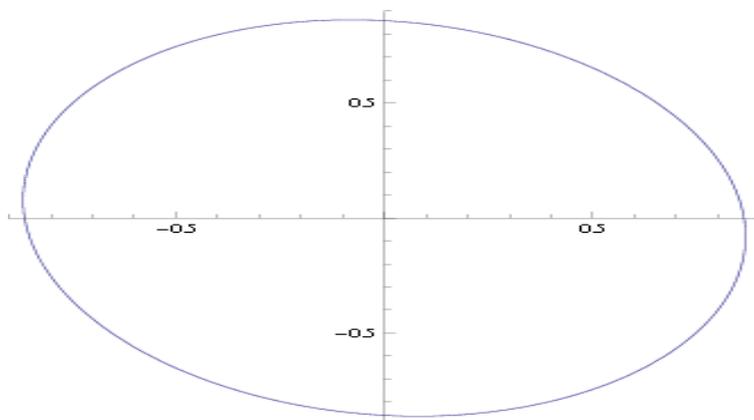}}
\caption{Profile curve for $m=2$, $H=0.57$, in this case $C=3.5313222039296357$}
\end{figure}

\vfill
\eject

\begin{figure}[h]\label{profiles m2}
\centerline{\includegraphics[width=10cm,height=5.5cm]{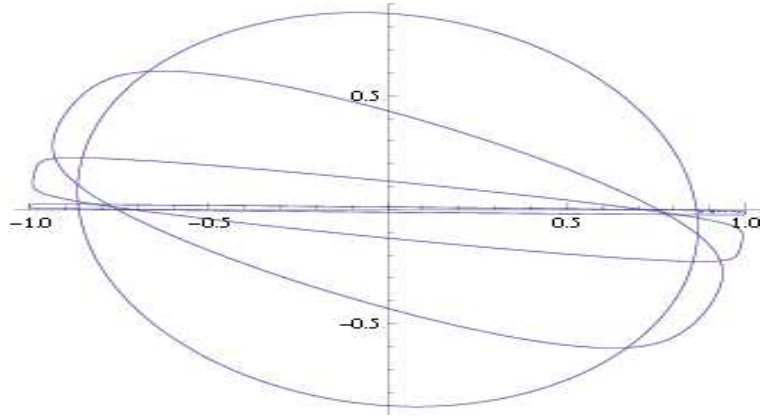}}
\caption{Profile curve for $m=2$, $H=0.001,\, H=0.1,\, H=0.3,\, H=0.57$.}
\end{figure}

\vskip2cm

\begin{figure}[h]\label{stereographic h=0.1}
\centerline{\includegraphics[width=10cm,height=5.5cm]{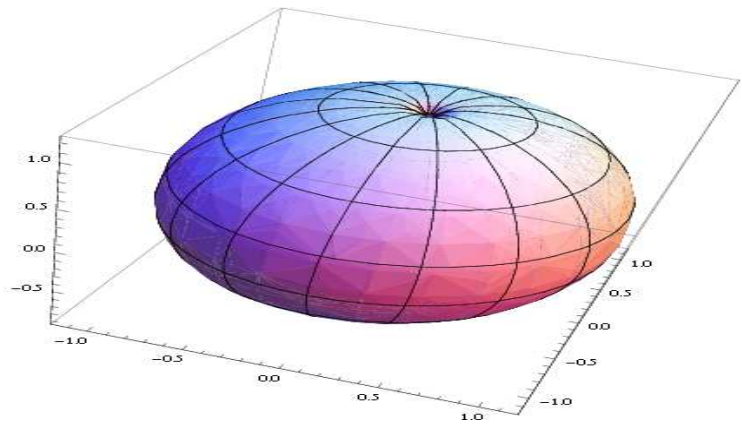}}
\caption{Stereographic projection for the surface with CMC $H=0.1$}
\end{figure}

\vfill
\eject

\begin{figure}[h]\label{stereographic half h=0.1}
\centerline{\includegraphics[width=10cm,height=5.5cm]{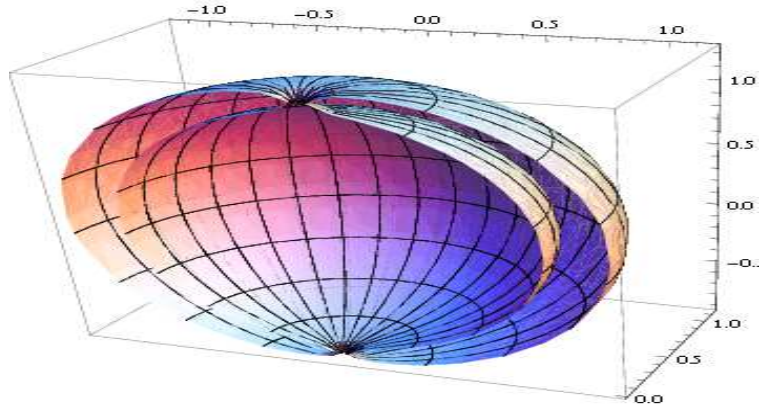}}
\caption{Stereographic projection of half the surface with CMC $H=0.1$}
\end{figure}

\vskip2cm

\begin{figure}[h]\label{neck h=0.1}
\centerline{\includegraphics[width=10cm,height=5.5cm]{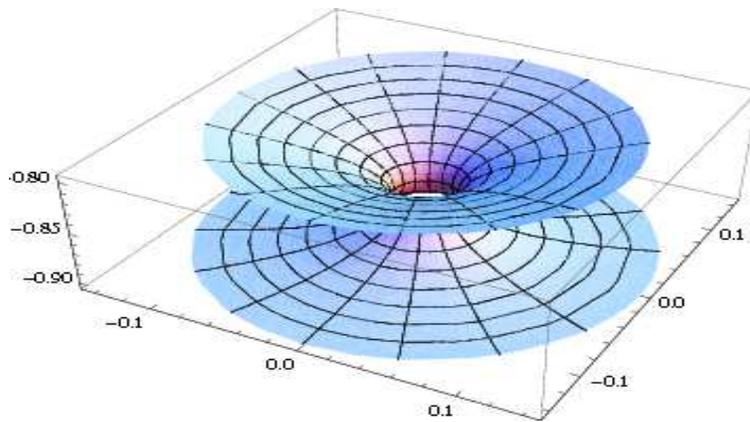}}
\caption{Stereographic projection one of the two catenoid necks of the surface with CMC $H=0.1$}
\end{figure}

\vfill
\eject

\begin{figure}[h]\label{stereographic h=0.3}
\centerline{\includegraphics[width=10cm,height=5.5cm]{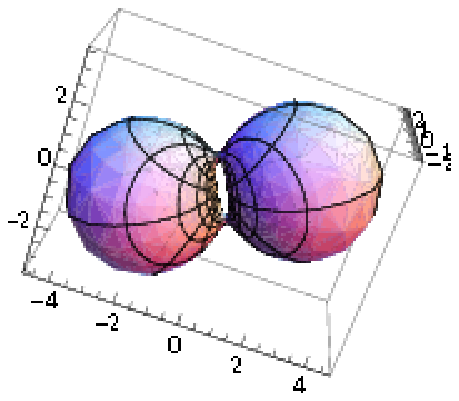}}
\caption{Stereographic projection of the surface with CMC $H=0.3$ and $m=2$}
\end{figure}

\vskip2cm

\begin{figure}[h]\label{stereographic h=0.57}
\centerline{\includegraphics[width=10cm,height=5.5cm]{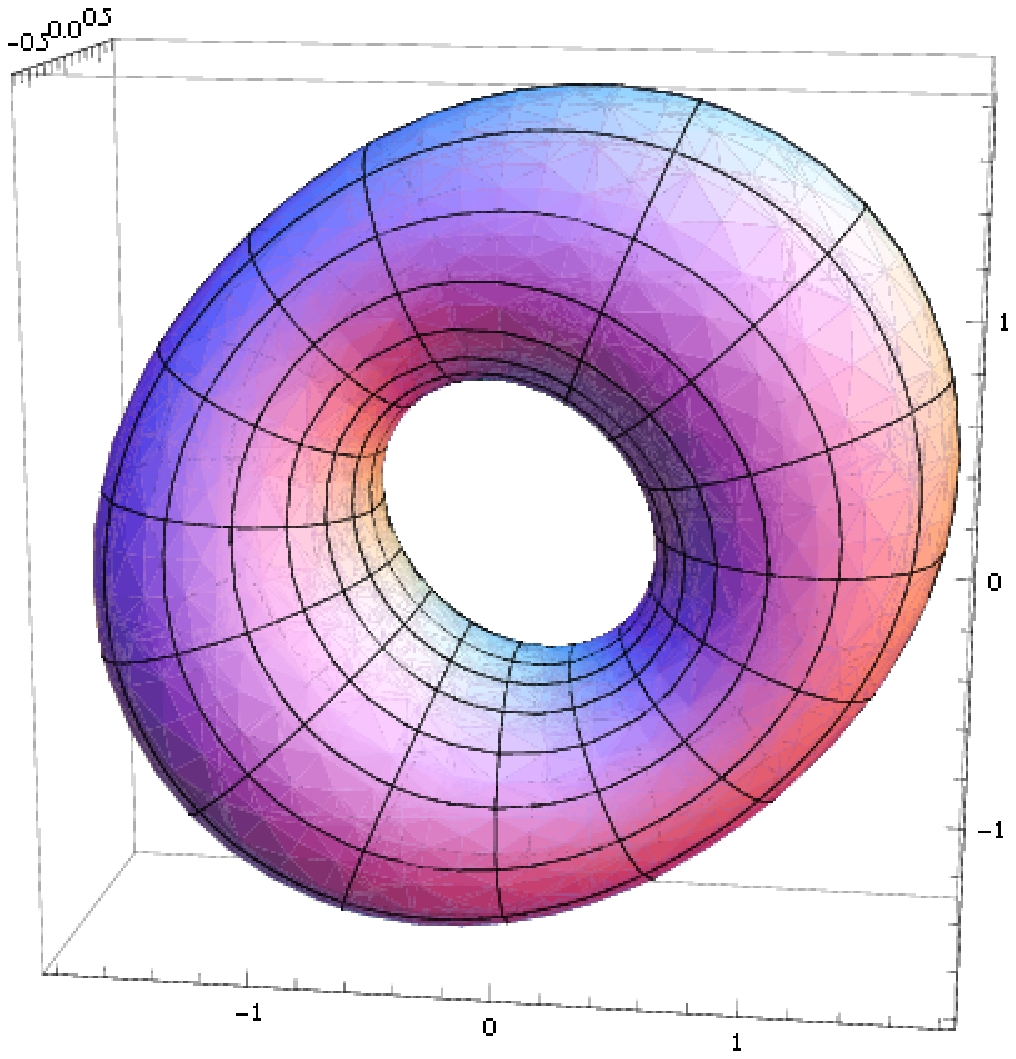}}
\caption{Stereographic projection of the surface with CMC $H=0.57$ and $m=2$}
\end{figure}

\vfill
\eject

\begin{figure}[h]\label{profile h=0.5774}
\centerline{\includegraphics{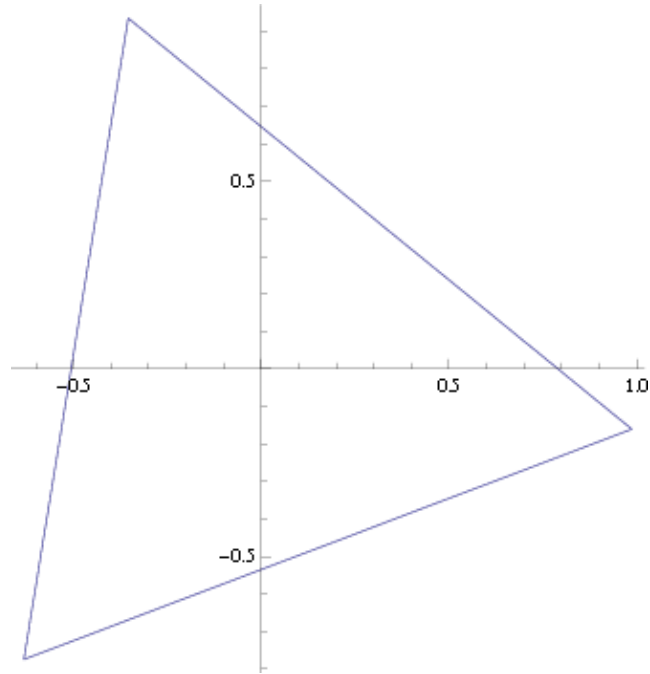}}
\caption{Profile curve for $m=3$ and $H=0.5774$, in this case $C=346879.6632142387$}
\end{figure}

\begin{figure}[h]\label{profile h=0.6}
\centerline{\includegraphics{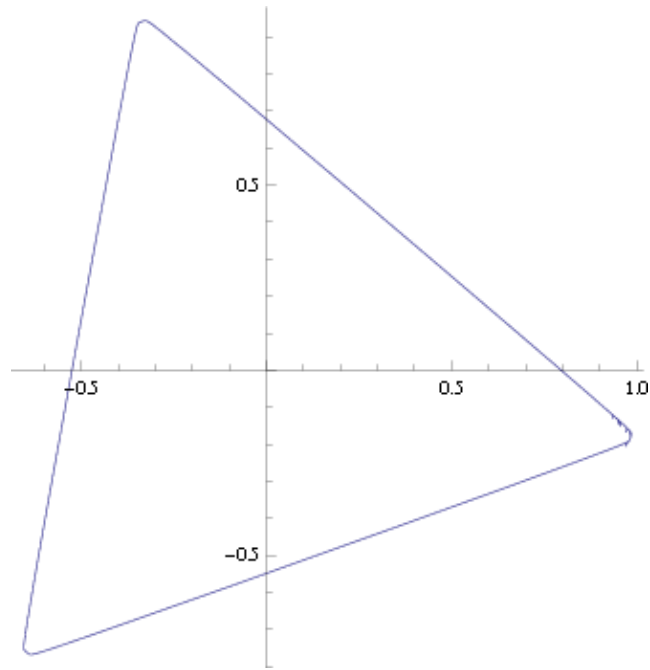}}
\caption{Profile curve for $m=3$ and $H=0.6$, in this case $C=365.3705636110441$}
\end{figure}

\vfill
\eject

\begin{figure}[h]\label{profile h=0.8}
\centerline{\includegraphics{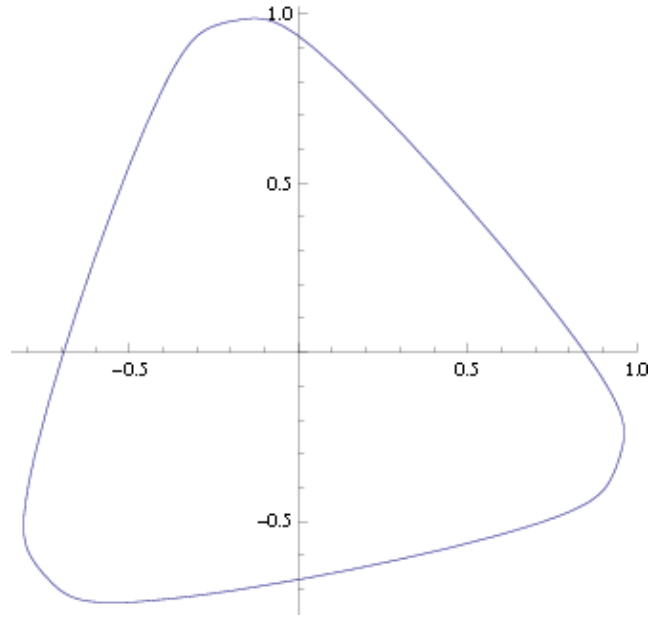}}
\caption{Profile curve for $m=3$ and $H=0.8$, in this case $C=22.320379289179478$}
\end{figure}

\begin{figure}[h]\label{profile h=1.0}
\centerline{\includegraphics{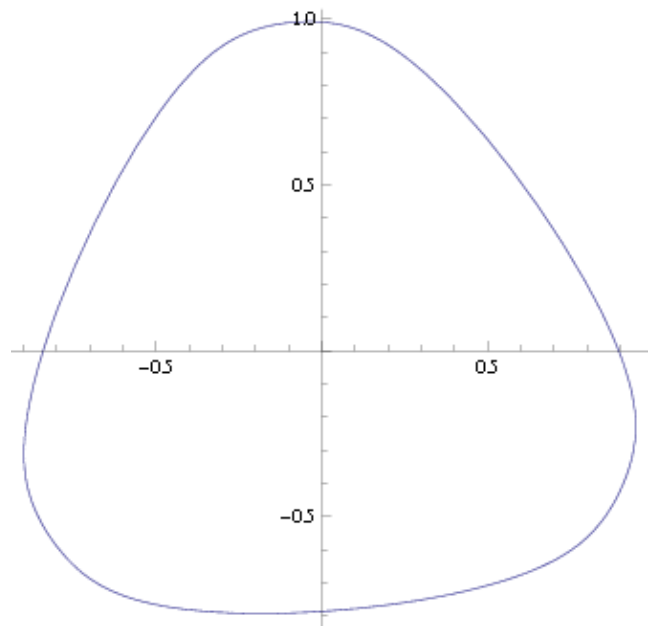}}
\caption{Profile curve for $m=3$ and $H=1.0$, in this case $C=9.908469426660892$}
\end{figure}

\vfill
\eject

\begin{figure}[h]\label{profile h=1.2}
\centerline{\includegraphics{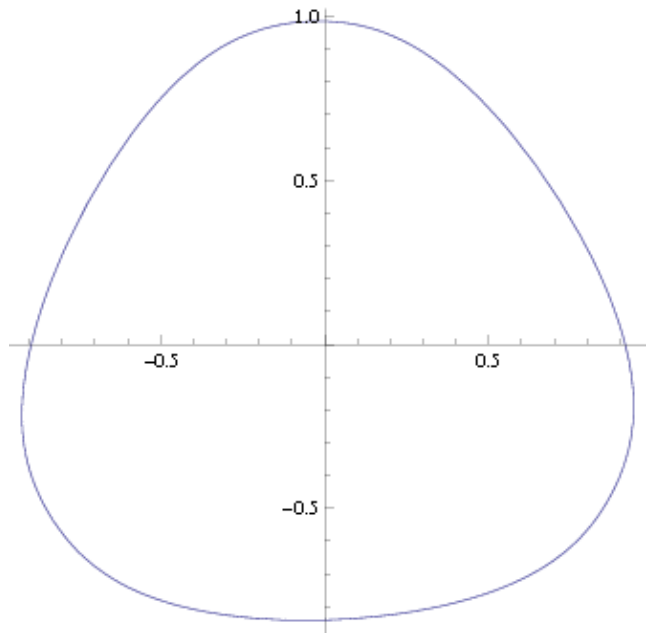}}
\caption{Profile curve for $m=3$ and $H=1.2$, in this case $C=6.084010495710457$}
\end{figure}

\begin{figure}[h]\label{profile h=1.237}
\centerline{\includegraphics{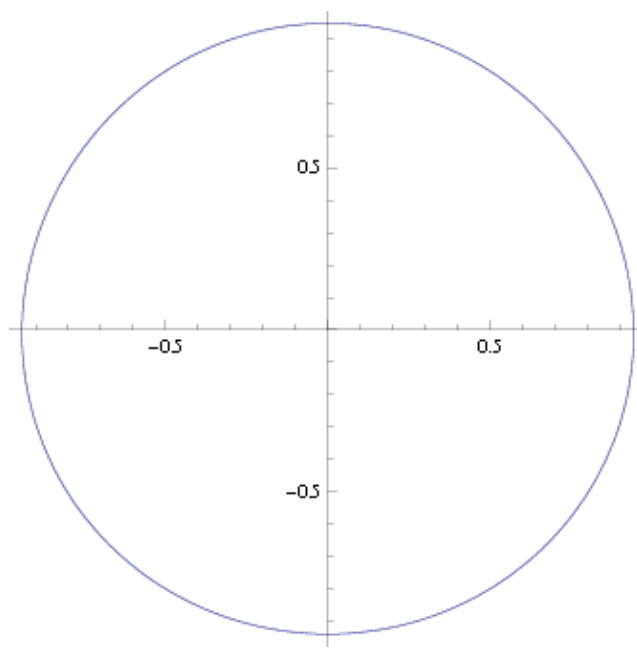}}
\caption{Profile curve for $m=3$ and $H=1.237$, in this case $C=5.6615177218839605$}
\end{figure}

\vfill
\eject

\begin{figure}[h]\label{profile curves m3}
\centerline{\includegraphics{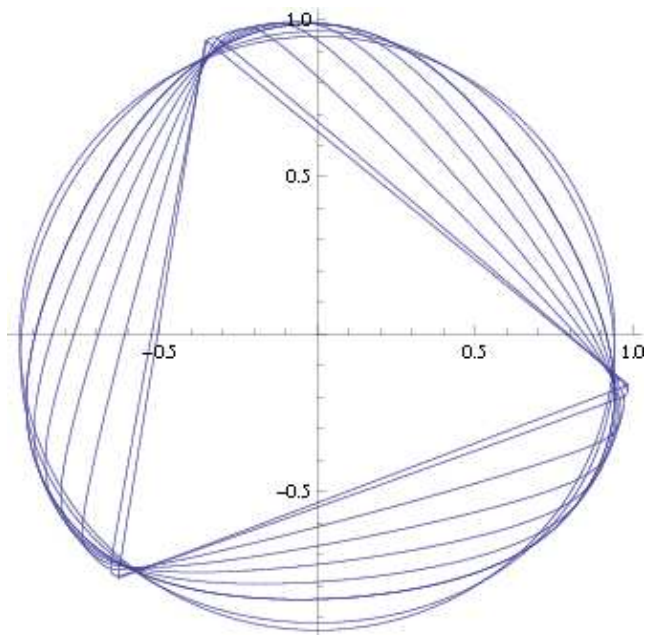}}
\caption{Profile curve for $m=3$, $H=0.5774,\, H=0.6,\, H=0.7,\, H=0.8,$
$H=1.0\, H=1.1\, H=1.2,\, H=1.22\, H=1.237$.}
\end{figure}

\vfill
\eject

\begin{figure}[h]\label{stereographic projection h=0.5774}
\centerline{\includegraphics{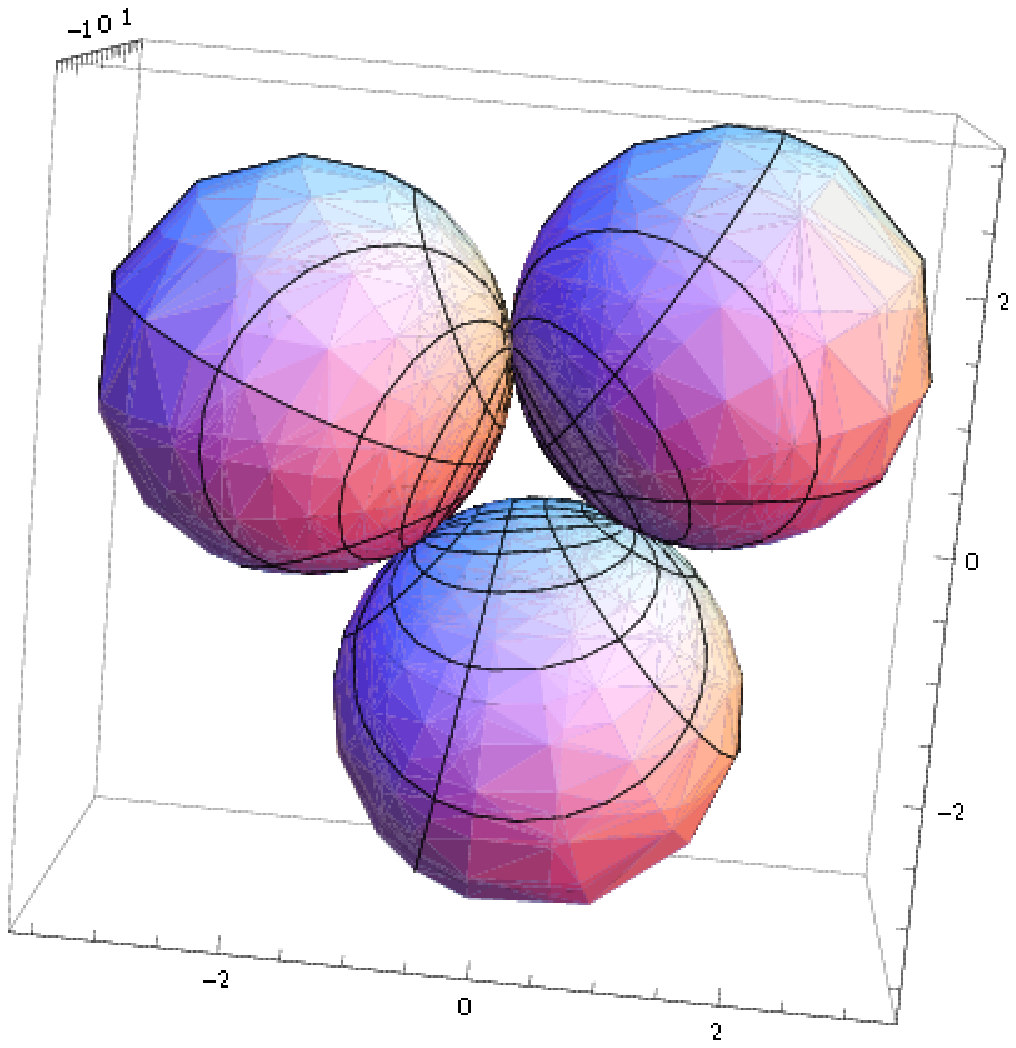}}
\caption{Stereographic projection of a surface with CMC $H=0.5774$ and $m=3$}
\end{figure}

\vfill
\eject

\begin{figure}[h]\label{mspm3h0p8}
\centerline{\includegraphics{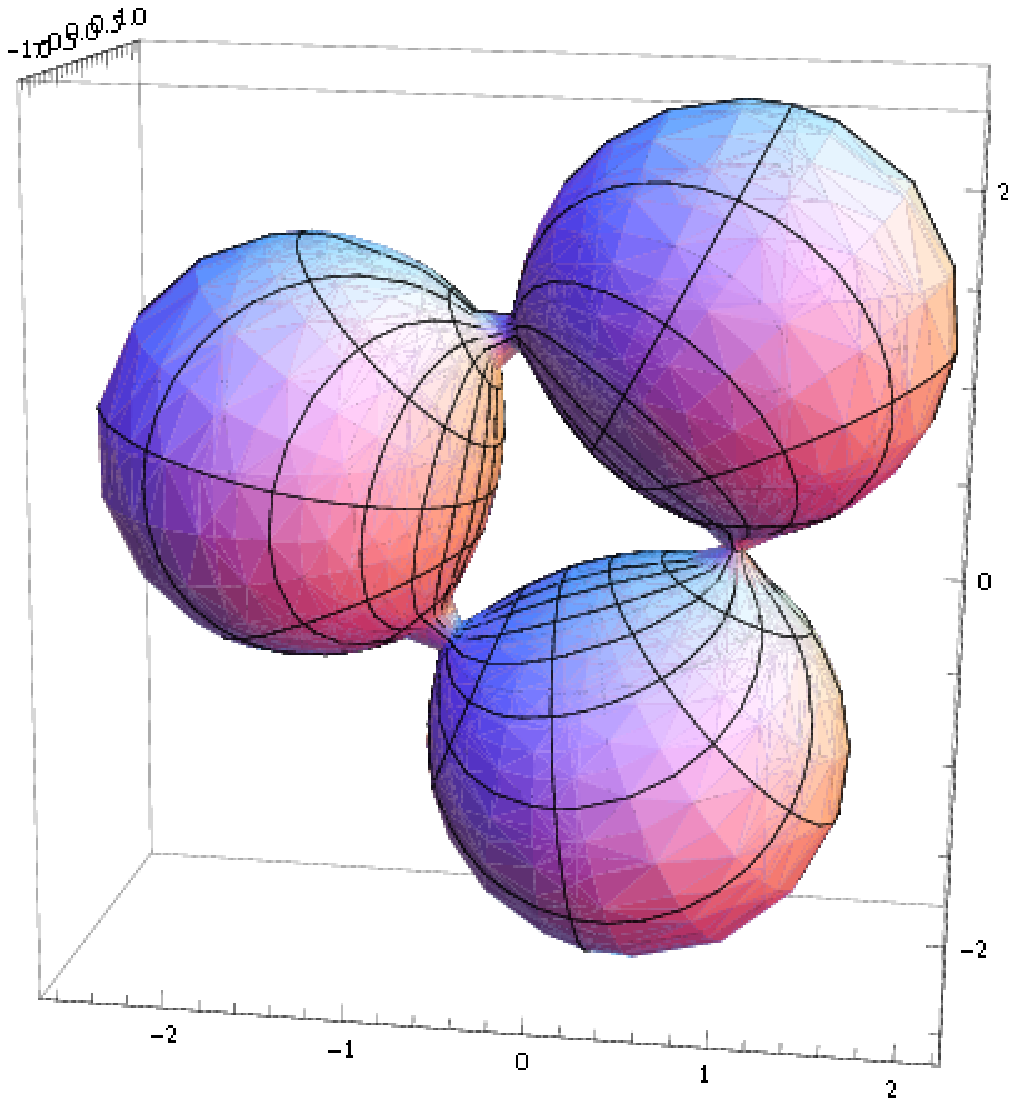}}
\caption{Stereographic projection of a surface with CMC $H=0.8$ and $m=3$}
\end{figure}

\vfill
\eject

\begin{figure}[h]
\centerline{\includegraphics{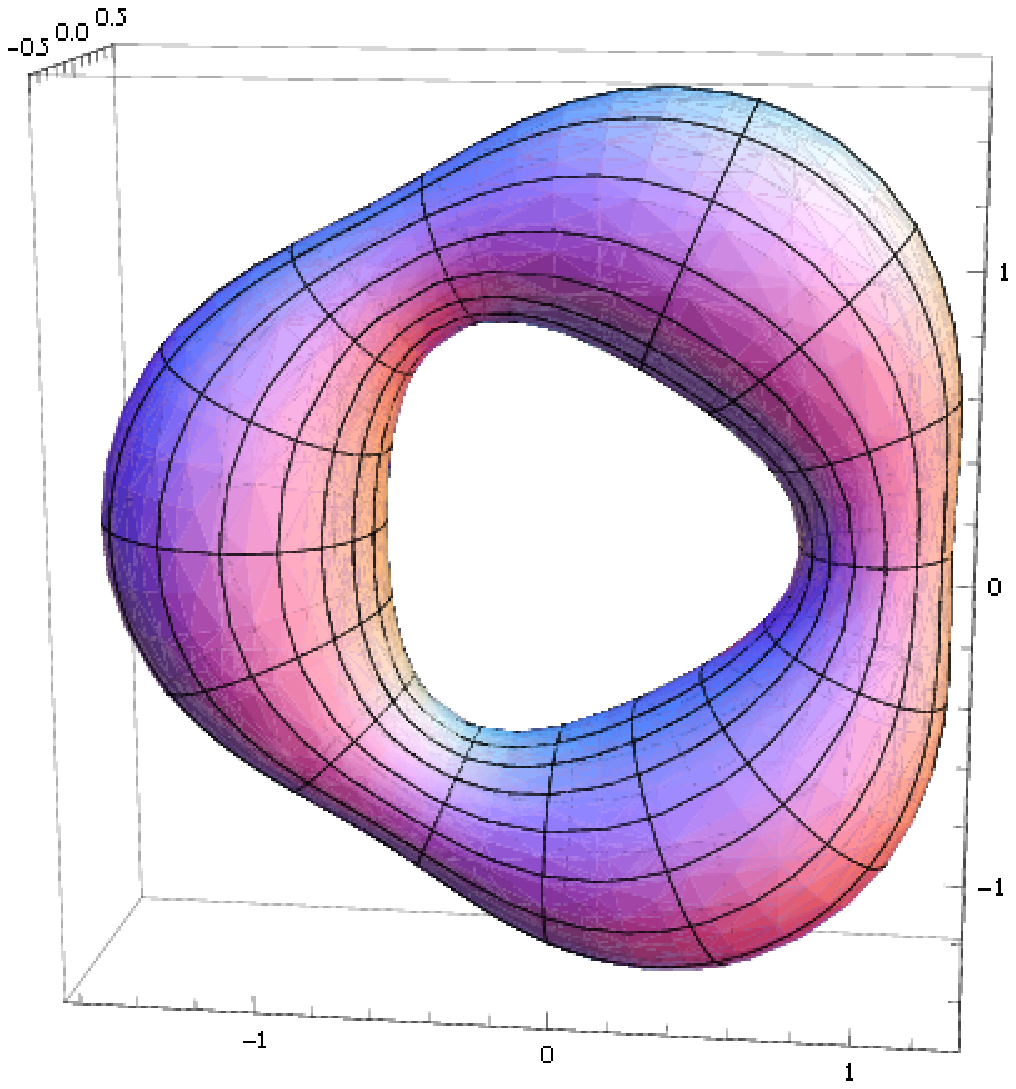}}
\caption{Stereographic projection of a surface with CMC $H=1.2$ and $m=3$}
\end{figure}

\vfill
\eject

\begin{figure}[h]
\centerline{\includegraphics{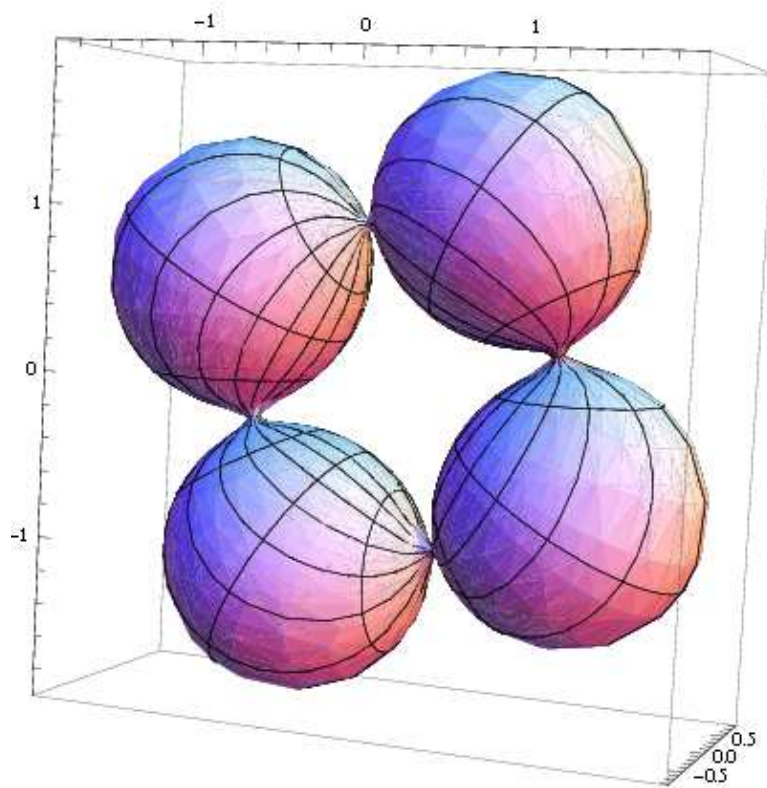}}
\caption{Stereographic projection a surface with CMC $H=1.2$ and $m=4$}
\end{figure}

\vfill
\eject

\subsection{Embedded solutions in hyperbolic spaces.}

In this section we will point out that the theorem above can be adapted to hyperbolic spaces. In this case we obtained the embedded hypersurfaces with not much effort since the Hyperbolic space is not compact. Here we will be considering the following model of the hyperbolic space,

$$H^{n+1}=\{\, x\in \bfR{{n+2}}: x_1^2+\dots+x_{n+1}^2-x_{n+2}^2=-1\, \} $$

The following notation will only be considered in this subsection. For any pair of vectors $v=(v_1,\dots,v_{n+2})$ and  $w=(w_1,\dots,w_{n+2})$, $\<v,w\>=v_1w_1+v_{n+1}w_{n+1}-v_{n+2}w_{n+2}$.

\begin{thm}

Let $g_{C,H}:{\bf R}\to {\bf R}$ be a positive solution of the equation

\begin{eqnarray}\label{eq hyperbolic}
 (g^\prime)^2+g^{2-2n}+(H^2-1)g^2+2Hg^{2-n}=C
 \end{eqnarray}

associated with a non negative $H$  and a positive constant $C$. If $\mu,\lambda,r,\theta:{\bf R}\to {\bf R}$  and are defined by

$$r=\frac{g_{C,H}}{\sqrt{C}},\quad  \lambda=H+g_{C,H}^{-n},\, \mu=nH-(n-1)\lambda=H-(n-1)g_{C,H}^{-n} \com{and}\theta(u)=\int_0^u\frac{r(s)\lambda(s)}{1+r^2(s)}ds$$

then, the map $\phi:S^{n-1}\times {\bf R}\to H^{n+1}$ given by

\begin{eqnarray}\label{the immersions hyperbolic}
\phi(y,u)=(\, r(u)\, y,\sqrt{1+r(u)^2}\, \sinh(\theta(u)),\sqrt{1+r(u)^2} \, \cosh(\theta(u))\, )
\end{eqnarray}

defines an embedded hypersurface in $H^{n+1}$ with constant mean curvature $H$. Moreover, if
$H^2>1$,  the embedded manifold defined by $(\ref{the immersions hyperbolic})$ admits the group $O(n)\times Z$ in its group of isometries, where $Z$ is the group of integers.
\end{thm}

{\bf Remark:}
{\sl Arguments similar to those in section (\ref{section the function w and its solutions}) show that  it is not difficult to find positive values $C$ that lead to positive  solutions of the equation (\ref{eq hyperbolic}) in terms of the inverse of a function defined by an integral.
}

\begin{proof}

A direct computation shows the following identities,

$$ (r^\prime)^2+\lambda^2=1+r^2,\com{and} \lambda r^\prime +r\lambda^\prime=\mu r^\prime$$

Let us define

$$B_2(u)=(0,\dots,0,\sinh(\theta(u)),\cosh(\theta(u)))\com{and} B_3(u)=(0,\dots,0,\cosh(\theta(u)),\sinh(\theta(u)))$$

Notice that $\<B_2,B_2\>=-1$,  $\<B_3,B_3\>=-1$, $\<B_2,B_3\>=0$,  $B_2^\prime=\frac{r\lambda}{1+r^2}B_3$ and
$B_3^\prime=\frac{r\lambda}{1+r^2}B_2$, moreover, we have that the map $\phi$ can be written as

$$\phi=r(y,0,0)+\sqrt{1+r^2}\, B_2$$

A direct verification shows that $\<\phi,\phi\>=-1$ and that

$$\frac{\partial{\phi}}{\partial u}=r^\prime\, (y,0,0)+\frac{rr^\prime}{\sqrt{1+r^2}} \, B_2+
\frac{r \lambda}{\sqrt{1+r^2}} \, B_3$$

is a unit vector, i.e, $\<\frac{\partial{\phi}}{\partial u},\frac{\partial{\phi}}{\partial u}\>=1$.
We have that the tangent space of the immersion at $(y,u)$ is given by

$$T_{\phi(y,u)}=\{(v,0,0)+s\, \frac{\partial \phi}{\partial u}: \<v,y\>=0\com{and} s\in{\bf R}\}$$

A direct verification shows that the map

$$\nu=-r\lambda\,   (y,0,0)-\frac{r^2\, \lambda}{\sqrt{1+r^2}} \, B_2 + \frac{r^\prime}{\sqrt{1+r^2}} \, B_3$$

satisfies that $\<\nu,\nu\>=1$, $\<\nu,\frac{\partial \phi}{\partial u}\>=0$ and for any $v\in\bfR{n}$ with $\<v,y\>=0$ we have that $\<\nu,(v,0,0)\>=0$. It then follows that $\nu$ is a Gauss map of the immersion $\phi$. The fact that the immersion $\phi$ has constant mean curvature $H$ follows because for any unit vector $v$ in $\bfR{n}$ perpendicular to $y$, we have that

$$\beta(t)=(r\cos(t)y+r\sin(t)v,0,0)+\sqrt{1+r^2}\, B_2=\phi(\cos(t)y+r\sin(t)v,u)$$

satisfies that $\beta(0)=\phi(y,u)$, $\beta^\prime(0)=rv$ and

$$\frac{d \nu(\beta(t))}{dt}\big{|}_{t=0} = d\nu(rv)=-r\lambda\, v$$

Therefore, the tangent vectors of the form $(v,0,0)$ are principal directions with principal curvature $\lambda$ and multiplicity $n-1$. Now, since $\<\frac{\partial \nu}{\partial u},(v,0,0)\>=0$, we have that $\frac{\partial \phi}{\partial u}$ defines a principal direction, i.e. we must have that $\frac{\partial \phi}{\partial u}$ is a multiple of $\frac{\partial \phi}{\partial u}$.  A direct verification shows that,

$$\<\frac{\partial \nu}{\partial u},y\>=-\lambda^\prime \, r-\lambda r^\prime=-\mu\, r^\prime=-(nH-(n-1)\lambda)r^\prime$$

We also have that $\<\frac{\partial \phi}{\partial u},y\>=r^\prime$, therefore,

$$\frac{\partial \nu}{\partial u}= d\nu(\frac{\partial \phi}{\partial u})=-\mu\, \frac{\partial \phi}{\partial u}=-(nH-(n-1)\lambda) \frac{\partial \phi}{\partial u}$$

It follows that the other principal curvature is $nH-(n-1)\lambda$. Therefore $\phi$ defines an immersion with constant mean curvature $H$, this proves the first item in the Theorem.
This immersion is  embedded because the immersion $\phi$ is one to one as we can easily check using the fact that  whenever $H\ge0$, the function $\theta$ is strictly increasing. In order to prove the condition on the isometries of the immersion when $H>1$ we notice first that the ODE (\ref{eq hyperbolic}) can be written as

$$(g^\prime)^2=g^{2-2n}\, q(g)\com{where} q(v)=Cv^{2n-2}-(H^2-1) v^{2n}-2H v^n-1$$

Since $q(0)=-1$ and the leading coefficient of $q$ is negative under the assumption that $H>1$, then by the arguments used in section (\ref{section the function w and its solutions}) we conclude that a positive solution $g$ of (\ref{eq hyperbolic}) must be periodic, moreover the values of $g$ must move from two positive roots $t_1$ and $t_2$. Now if $T$ is the period of $g$ and we define

$$ K=\int_0^T\frac{r(u)\lambda(u)}{1+r^2(u)}du$$

then we have,

$$ \com{For any integer $j$ and $u\in [jT,(j+1)T]$ we have that }\theta(u)=jK+\theta(u-jT) $$

Using the equation above we get that the immersion $\phi$ is invariant under the group generated by hyperbolic rotations of the angle $K$ in the $x_{n+1}$-$x_{n+2}$ plane. This concludes the theorem.

\end{proof}

\subsection{Solutions in Euclidean spaces.}

In this section we will point out that the same kind of theorem can be adapted to Euclidean spaces. In this case we obtain the embedded and not embedded Delaunay hypersurfaces.

\begin{thm}

Let $g_{C,H}:{\bf R}\to {\bf R}$ be a positive solution of the equation

\begin{eqnarray}\label{eq hyperbolic}
 (g^\prime)^2+g^{2-2n}+H^2\, g^2+2Hg^{2-n}=C
 \end{eqnarray}

associated with a real number  $H$  and a positive constant $C$. If $\mu,\lambda,r,R:{\bf R}\to {\bf R}$  and are defined by

$$r=\frac{g_{C,H}}{\sqrt{C}},\quad  \lambda=H+g_{C,H}^{-n},\, \mu=nH-(n-1)\lambda=H-(n-1)g_{C,H}^{-n} \com{and}R(u)=\int_0^u r(s)\lambda(s)ds$$

then, the map $\phi:S^{n-1}\times {\bf R}\to \bfR{{n+1}}$ given by

\begin{eqnarray}\label{the immersions euclidean}
\phi(y,u)=(\, r(u)\, y,R(u) )
\end{eqnarray}

defines an immersed hypersurface in $\bfR{{n+1}}$ with constant mean curvature $H$. Moreover, if
$H\ge 0$,  the manifold defined by $(\ref{the immersions euclidean})$ is embedded. We also have that when $n>2$, up to rigid motions they are the only CMC hypersurfaces with exactly two principal curvatures.
\end{thm}

\begin{proof}

A direct computation shows the following identities,

$$ (r^\prime)^2+\lambda^2=r^2,\com{and} \lambda r^\prime +r\lambda^\prime=\mu r^\prime$$

In this case we have that the map

$$\nu(y,u)=(-r(u)\lambda(u)\,  y,r^\prime(u))$$

is a Gauss map of the immersion. A direct computation shows that indeed this immersion has constant mean curvature $H$. The fact that the immersion is an embedding when $H\ge0$ follows from the fact that $\lambda>0$ in this case and therefore the function $R$ is strictly increasing. For the last part of the theorem we will use the same notation used in the previous sections, and in particular we define the functions $w,\lambda$ on the whole manifold as before, and we extend the function $r$ to the manifold by defining it as $r=\frac{w}{\sqrt{c}}$. We have that,

\begin{enumerate}
\item
the vector $\lambda\, r e_n+ e_n(r) \, \nu$ is a unit constant vector on the whole manifold, we can assume that this vector is the vector $(0,\dots,0,1)$
\item
The vector $\eta=-e_n(r)\, e_n +\lambda\, r\, \nu$ is constant along the geodesic defined by the vector field $e_n$, i.e, we can prove that $\bar{\nabla}_{e_n}\eta$ vanishes.
\item
From the last items we can solve for $e_n$ in terms of the vectors $(0,\dots,0,1)$ and $\eta$, and then, integrate in order to get the profile curves.
\item

The vector field $x+r\eta$ is independent of the integral submanifolds of the distribution $\hbox{Span}\{e_1,\dots,e_{n-1}\}$.

\item

The previous considerations and the fact that the vectors $e_1,\dots e_{n-1}$ are perpendicular to the vector $\eta$ and $(0,\dots,0,1)$ imply that the integral submanifolds of the distribution $\hbox{Span}\{e_1,\dots,e_{n-1}\}$ are spheres with center at $x+r\eta$ and radius $r$. Notice that $||x-(x+r\eta)||=r$.

\item

If we fix a point $p_0$ and we define the geodesic $\gamma(u)$ as before, then, without loss of generality we may assume that $\eta(p_0)=(0,\dots,0,1,0)=\eta(u)$ and therefore, we can also assume by doing a translation, if necessary, that

$$\gamma(u)=\int_0^ue_n(u)=\int_0^u(0,\dots,0,-r^\prime(u),\lambda(u)\, r(u))=
(0,\dots,-r(u),R(u))$$

Where $R(u)=\int_0^u\lambda(t)\, r(t) dt$. The theorem follows by noticing that the center of the integral submanifolds take the form $\gamma(u)+r(u)\eta(\gamma(u))=(0,\dots,0,R(u))$
\end{enumerate}

\end{proof}

 In the case $n=2$ we can find explicit solutions. For any positive $C>4H$, they look like,

$$\phi(u,v)=(\,r(u)\cos(v),r(u)\sin(v),R(u)) $$

where ,

$$R(u)=\int_0^u \frac{C + \sqrt{C (C - 4 H)} \cos(2 H y)\, }{\sqrt{2C} \sqrt{
 C - 2 H + \sqrt{C (C - 4 H)} \cos(2 H y)\,}} \, dy\com{and}r(u)=\frac{\sqrt{C - 2 H + \sqrt{C (C - 4 H)} \cos(2 H u)}}{\sqrt{2} \sqrt{C} H}$$

Here there is the graph of a non embedded Delaunay surface,

\begin{figure}[h]
\centerline{\includegraphics{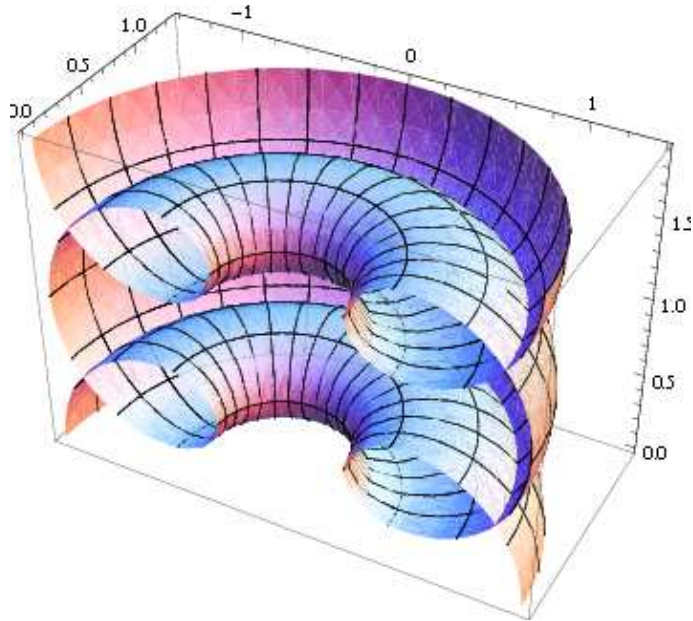}}
\caption{Half rotation of a non embedded Delaunay surface with CMC $H=-1$, here $C=2$}
\end{figure}

\vfill
\eject

\end{document}